\documentclass[12pt]{article}
\usepackage{graphics,graphicx}
     \usepackage{amsmath}

\usepackage{amssymb}
\usepackage[latin1]{inputenc}
\usepackage{mathrsfs}

\def\fullpage {
\addtolength{\topmargin}{-1.5 cm}
\addtolength{\oddsidemargin}{-1.5 cm}
\addtolength{\textwidth}{+3 cm}
\addtolength{\textheight}{+3 cm} }

  \fullpage

        \newcommand{\lab}[1]{\label{#1}}

\catcode`@=11
\@addtoreset{equation}{section}

\catcode`@=12

\def\blackslug{\hbox{\kern1pt\vrule height6pt width4pt  depth1pt\kern1pt}}
\def\qed{\penalty 500\hbox{\quad\blackslug}\ifmmode\else\par
                                             \vskip4.5pt plus3pt minus2pt\fi}
\newtheorem{thm}{Theorem}

\newtheorem{lemma}[thm]{Lemma}

\newtheorem{prop}[thm]{Proposition}

\def\proof{\par\noindent{\bf Proof.\enspace}\rm}
\newcommand{\proofof}[1]{\par\noindent{\bf Proof of #1.\enspace}\rm}

   \def\G{{\cal G}}
   \def\E{{\cal E}}
   \def\P{{\cal P}}

\def\eps{\epsilon}

\def\ex{{\bf E}}
\def\pr{{\bf P}}

\def\BZ{{\mathbb Z}}
\def\BC{{\mathbb C}}

\def\Gnd{\G_{n,d}}
\def\Gnfour{\G_{n,4}}

\def\Gnsix{\G_{n,6}}
\def\Gnten{\G_{n,10}}
\def\Pnd{\P_{n,d}}
\def\Pn5{\P_{n,5}}
\def\Gnp{{\cal G}(n,p)}

\def\st{\; | \:}

\def\be{\begin{equation}}
\def\ee{\end{equation}}
\def\bea{\begin{eqnarray}}
\def\eea{\end{eqnarray}}
\def\bean{\begin{eqnarray*}}
\def\eean{\end{eqnarray*}}

\newcommand{\eqn}[1]{(\ref{#1})}
\newcommand{\bel}[1]{\be\lab{#1}}

\def\transpose{\top}
\def\mvec{{\bf vec}}
\def\mdet{{\rm det}}
\def\onesvec{{\mathbf{1}}}
\def\vtheta{{\bf \theta}}
\def\Z{\BZ}
\def\M{{\cal M}}

\newcommand{\real}{\ensuremath {\mathbb R} }

\newcommand{\mbf}[1] {\text{\boldmath$#1$}}
\newcommand{\remove}[1] {}

\DeclareMathOperator{\poly}{poly}

\newcommand{\cB} {\ensuremath{\mathcal B}}
\newcommand{\cG} {\ensuremath{\mathcal G}}
\newcommand{\cL} {\ensuremath{\mathcal L}}
\newcommand{\cM} {\ensuremath{\mathcal M}}
\newcommand{\cP} {\ensuremath{\mathcal P}}

\newcommand{\hcL} {\ensuremath{\widehat{\mathcal L}}}

\newcommand{\B} {\ensuremath{\mathsf B}}

\begin{document}
  \title{ On the chromatic number of random $d$-regular graphs}
\author{
Graeme~Kemkes\thanks{Partially supported by NSERC CGS and PDF awards. Some of this research appeared in this author's PhD thesis at the University of Waterloo.}\\
{\small Department of Mathematics}\\
{\small University of California, San Diego}\\
{\small La Jolla, CA 92093-0112, US}\\
{\small {\tt gdkemkes@uwaterloo.ca }}
\and
Xavier~P\'erez-Gim\'enez\thanks{Partially supported by the Province of Ontario under the Post-Doctoral Fellowship (PDF) Program.}\\
{\small Department of Combinatorics and Optimization} \\
{\small University of Waterloo}\\ 
{\small Waterloo ON, Canada}\\
{\small {\tt xperez@uwaterloo.ca }}
\and
Nicholas Wormald\thanks{Supported by the  Canada Research Chairs Program and NSERC.}\\
{\small Department of Combinatorics and Optimization} \\
{\small University of Waterloo}\\ 
{\small Waterloo ON, Canada}\\
{\small {\tt nwormald@uwaterloo.ca }} 
} 
\date{}
\maketitle
\begin{abstract}
In this work we show that, for any fixed $d$, random $d$-regular graphs asymptotically almost surely can be coloured with $k$ colours, where $k$ is the smallest integer satisfying $d < 2(k-1)\log(k-1)$. From previous lower bounds due to Molloy and Reed, this establishes the chromatic number to be asymptotically almost surely $k-1$ or $k$. If moreover $d>(2k-3)\log(k-1)$, then the value $k-1$ is discarded and thus the chromatic number is exactly determined. Hence we improve a recently announced result by Achlioptas and Moore in which the chromatic number was allowed to take the value $k+1$. Our proof applies the small subgraph conditioning method to the number of balanced $k$-colourings, where a colouring is  \emph{balanced} if the number of vertices of each colour is equal.
\end{abstract}
%
\section{Introduction}
%
The chromatic number $\chi$ of random graphs is a topic that has attracted considerable interest since the breakthrough achieved by Shamir and Spencer~\cite{SS}, which marked one of the first applications of martingales in combinatorics.
For the classical Erd\H os-R\' enyi model $\Gnp$, a celebrated result by Bollob{\'a}s~\cite{boll88} later extended by {\L}uczak~\cite{L91} showed that if $pn\to\infty$ then asymptotically almost surely (a.a.s.)
\[
\chi(\Gnp)\sim \frac{n \log(1/(1-p))}{2\log(np)}.
\]
Here and in similar statements, an event occurs a.a.s.\ if its probability tends to 1 as $n$ tends to infinity.
For $p=c/n$, Achlioptas and Naor~\cite{AN} proved that the chromatic number of $\Gnp$ is a.a.s.\ $k$ or $k+1$ where $k$ is the smallest positive integer with $2k \log k > c$. Moreover, they discarded the case $k$ for roughly half of the values of $c$.
In the same direction, Coja-Oghlan, Panagiotou and Steger~\cite{CPS} showed that a.a.s.\ $\chi(\Gnp)\in\{ k,k+1,k+2\}$ for $p<n^{-3/4 -\eps}$    
where $k$ is the smallest positive integer satisfying $2k \log k > p(n-1)$.
Meanwhile, some other results gave concentration of the chromatic number without determining the values so precisely:
{\L}uczak~\cite{L91b} proved that $\chi(\Gnp)$ is a.a.s.\ 2-point-concentrated if $p<n^{-5/6-\eps}$, and later Alon and Krivelevich~\cite{AK97} extended this to $p<n^{-1/2-\eps}$.

More recently, results have been published about the chromatic number for the model
$\Gnd$ of random $d$-regular graphs, which is the probability space on $d$-regular graphs with $n$ vertices having uniform distribution. For basic results and notation on random regular graphs, see~\cite{NW99}. Hereinafter, $dn$ is always assumed to be even for feasibility.
For fixed $d$, Molloy and Reed~\cite{MR} showed that if $q(1-1/q)^{d/2}<1$ then $\chi(\Gnd)>q$ a.a.s.
Then, for $d<n^{1/3-\eps}$, Frieze and {\L}uczak~\cite{FL92} established that
\[
\chi(\Gnd) = \frac{d}{2\log d}+O\left(\frac{d \log \log d}{\log^2 d}\right),
\]
and later Cooper, Frieze, Reed and Riordan~\cite{CFRR02} extended the same asymptotic formula to apply to $d\le n^{1-\eps}$. Similarly, the range $n^{6/7+\eps}\le d\le 0.9n$ was covered by Krivelevich, Sudakov, Vu and Wormald~\cite{KSVW01}, who showed that $\chi(\Gnd) \sim \frac{n}{2\log_b d}$  a.a.s.\ where $b=n/(n-d)$.
Achlioptas and Moore \cite{AM} recently announced a significant new result for constant $d$.
They stated that if $k$ is the smallest integer satisfying $d < 2(k-1)\log(k-1)$
then a.a.s.\  $\chi(\Gnd)$ is 
$k-1$, $k$, or $k+1$. If, in addition, $d > (2k-3)\log(k-1)$, then a.a.s.\ $\chi(\Gnd)$ is $k$ or $k+1$. 
Finally, Ben-Shimon and Krivelevich~\cite{BK08} established 2-point concentration of $\chi(\Gnd)$ for $d=o(n^{1/5})$.

In this paper we restrict the set of possible values for the chromatic number given by Achlioptas and Moore, and show that $\chi(\Gnd)$ a.a.s.\ cannot be $k+1$. Therefore this reduces the range of possibilities for $\chi(\Gnd)$ to only a.a.s.\ 
$k-1$ and $k$, in the first case, and establishes that $\chi(\Gnd)=k$ a.a.s.\
in the second case.
In particular it provides
an alternate proof of the results of Shi and Wormald \cite{SW4, SW} 
that a.a.s.\ $\chi(\Gnfour)=3$  and $\chi(\Gnsix)=4$. It also establishes,
for example,
the previously-unknown result that a.a.s.\ $\chi(\Gnten)=5$ and $\chi(10^6)=46523$.
We essentially need to show that $\Gnd$ is a.a.s.\ $k$-colourable, since the abovementioned lower bound of Molloy and Reed implies that $\Gnd$ is a.a.s.\ not $(k-2)$-colourable, and for the second case, not $(k-1)$-colourable. Our  basic approach for the upper bound is similar to that of Achlioptas and Moore, in that we analyse the second moment of the number $Y$ of balanced $k$-colourings of random regular graphs. 
In fact, Achlioptas and Moore found that the second central moment $\ex (Y-\ex Y)^2$ is (essentially, in a more or less equivalent model of random graphs) asymptotically a non-zero constant times the first. Consequently, Chebyshev's inequality fails to show the result which we claim above. In cases like this, this failure of the second moment inequality
to establish $Y>0$ a.a.s.\
can, at least for random structures similar to $\Gnd$, be overcome by using 
the small subgraph conditioning method
of Robinson and the third author. (See~ \cite[Chapter 9]{JLR} 
and~\cite{NW99} for a full exposition of the method.)
Using this, we show  that $\Gnd$ is a.a.s.\ $k$-colourable.
\begin{thm}\label{thm:main}
Given any integer $d\ge3$, let $k$ be the smallest integer such that $d<2(k-1)\log(k-1)$. Then the chromatic number of $\cG_{n,d}$ is a.a.s.\ $k-1$ or $k$. If furthermore $d>(2k-3)\log(k-1)$, then the chromatic number of $\cG_{n,d}$ is a.a.s.\ $k$.
\end{thm}
Actually, almost all previous applications of the small subgraph conditioning method were for a random variable that counted large subgraphs in the random graph.
To apply the  method in
this setting we need to calculate the first and second moments of the number of balanced $k$-colourings ($n$ is then required to be divisible by $k$), 
as well as joint moments of the number of such colourings and the number 
of short cycles.
These computations are done in the well-known \emph{pairing} or {\em configuration model} $\Pnd$ which was first introduced
by  Bollob{\'a}s~\cite{Boll}. A \emph{pairing} in $\Pnd$
is a perfect matching on a set of $dn$ points which are grouped into $n$ cells of $d$ points each. A random pairing naturally corresponds  in an
obvious way to a random
$d$-regular multigraph (possibly containing loops or multiple edges),
in which each cell
becomes a vertex. Colourings of the multigraph then correspond to assignments of colours to the cells of the model.  The reader should refer to~\cite{NW99} for
aspects of the pairing model not explained here.
\begin{prop}
\label{prop:moments}
Fix integers $d,k\ge3$. Let $Y$ be the number of balanced $k$-colourings of a random $d$-regular
multigraph $\Pnd$ (where $n$ is restricted to the set of multiples of $k$).
\begin{description}
\item{(a)} For $m \geq 1$, let $X_m$ be the number of $m$-cycles
in $\Pnd$. Then
$$
\label{E1}
\ex Y \sim 
k^{k/2}
\left( \frac{k-1}{2\pi (k-2)} \right)^{(k-1)/2}
n^{-(k-1)/2}
k^n
\left( 1 - \frac{1}{k} \right)^{dn/2}
$$
and
\bel{E3}
\ex(Y[X_1]_{p_1}\cdots [X_j]_{p_j})
    \sim \prod_{m=1}^j \left( \lambda_m \left( 1 + \delta_m
\right)\right)^{p_m} \ex (Y)
\ee
where
\[
\lambda_m = \frac{(d-1)^m}{2m}  \qquad\mbox{and}\qquad 
\delta_m = \frac{ (-1)^m }{ (k-1)^{m-1} }.
\]
\item{(b)}
If furthermore $d<2(k-1)\log(k-1)$ then
$$
\label{E2}
\ex(Y^2) \sim 
\frac{k^k (k-1)^{k(k-1)}}{(k^2 -2k -d +2)^{(k-1)^2/2} (2\pi (k-2))^{k-1}}
n^{-(k-1)}
k^{2n}
\left( 1 - \frac{1}{k} \right)^{dn}.
$$
\end{description}
\end{prop}
We next compute
\bel{eSumld}
\sum_{m\geq 1} \lambda_m \delta_m^2 =
(k-1)^2
\log \left( \frac{ k-1 }{ \sqrt{ k^2 -2k -d +2 } } \right),
\ee
and verify that for $n$ divisible by $k$
\bel{eCondition}
\frac{\ex (Y^2)}{(\ex Y)^2} \sim
    \left( \frac{ k-1 }{ \sqrt{ k^2 -2k -d +2} }\right)^{ (k-1)^2 }
= \exp\left(\sum_{m\geq 1} \lambda_m \delta_m^2\right),
\ee
which is the last ingredient required for the application of the small subgraph conditioning method.
\proofof{Theorem~\ref{thm:main} (for $\mbf n$ divisible by $\mbf k$)}
Assume throughout the proof that $k$ divides $n$, and observe that all the conditions of Theorem~4.1 in~\cite{NW99} are verified by Proposition~\ref{prop:moments}, \eqref{eSumld} and~\eqref{eCondition}.
Thus we may apply the small subgraph conditioning method to conclude that
$\pr(Y>0\mid\E)\to1$,
where $\E=\bigwedge_{\delta_k=-1} \{X_k=0\} = \{X_1=0\}$
is the event of having no loops. Because $\pr\left( X_2 = 0 \right)$ is bounded away from $0$ for large $n$ (see~e.g.~\cite{NW99}), it follows that $Y > 0$ a.a.s.\ for the simple graphs in $\cG_{n,d}$, thus proving the required upper bound on the chromatic number. Now observe that from our choice of $k$ we have $d\ge 2(k-2)\log(k-2)>(2k-5)\log(k-2)$. Then the required lower bounds follow immediately from the fact that if $d>(2q-1)\log q$ then $\chi(\Gnd)>q$ (applied to $q=k-2$ or $q=k-1$ for each case in the statement). This is just a slightly weaker formulation of the result given by Molloy and Reed~\cite{MR}. (Their proof is reported in~\cite{SW}, Theorem~1.3.)\qed

The following two sections supply the proof of Proposition~\ref{prop:moments}. Finally, in
Section~\ref{sec:general_n}, we conclude the proof of Theorem~\ref{thm:main} by extending the argument to general $n$.

\section{Joint moments: proof of Proposition~\ref{prop:moments}$\mbf{(a)}$}
\label{sec:joint}

Let $Y$ be the number of balanced $k$-colourings of a random $d$-regular
multigraph $\Pnd$. For $m \geq 1$, let $X_m$ be the number of $m$-cycles
in $\Pnd$.
We estimate the expected value of $Y$ by enumerating all balanced
$k$-colourings of all multigraphs in $\Pnd$. There are 
$\binom{n}{n/k, n/k, \ldots, n/k}$ ways
to choose the $k$ colour classes. These choices are all equivalent so 
fix one. Suppose there are $b_{ij}=b_{ji}$ edges between
colour class $i$ and colour class $j$ (for $1\le i,j\le k$ and $i\ne j$). The colours of the
neighbours of all of the points of colour class $i$
can be then chosen in
$(dn/k)!/\prod_{1\le j \le k\atop j\ne i} b_{ij}!$ ways. After this determination
is made, edges are constructed by putting a perfect matching between
the corresponding points in each pair of classes, in one of 
$\prod_{1\le i<j \le k}b_{ij}!$
ways. Thus we have
\begin{eqnarray}
\left| \Pnd \right|
\ex (Y) &=& \binom{n}{n/k, n/k, \ldots, n/k}
\sum_{\{b_{ij}\}}
\left( \prod_{i=1}^k \frac{(dn/k)!}{\prod_{1\le j \le k\atop j\ne i} b_{ij}!} \right)
\prod_{1\le i<j \le k} b_{ij}! \nonumber \\
&=&
\binom{n}{n/k, n/k, \ldots, n/k}
(dn/k)!^k
\sum_{\{b_{ij}\}}
\frac{1}{\prod_{1\le i<j \le k} b_{ij}!} \nonumber \\
&=&
\binom{n}{n/k, n/k, \ldots, n/k}
(dn/k)!^k
\Big[\prod_{l=1}^k x_l^{dn/k}\Big] \prod_{1\le i<j \le k}\sum_{l\geq 0} \frac{(x_i x_j)^l}{l!} \nonumber \\
&=&
\binom{n}{n/k, n/k, \ldots, n/k}
(dn/k)!^k
\Big[\prod_{l=1}^k x_l^{dn/k}\Big] \exp\bigg(\sum_{1\le i<j \le k} x_i x_j\bigg),\nonumber
\end{eqnarray}
where square brackets denote the extraction of a coefficient from a generating function. A particular case of the following result gives us an accurate estimate of that coefficient. The proof is based on the saddlepoint method and is included later in this section.
%
%
%
%
\begin{lemma}
\label{lem:C(s)}
Let $k, d, a_1, a_2, \ldots, a_k$ be fixed integers with $k\geq 3$, $d$ positive,
and $s = \sum_{j=1}^k a_j$ even.
Let $C_{a_1, a_2,\ldots, a_k}$ denote the coefficient of
$x_1^{dn/k+a_1} x_2^{dn/k+a_2} \cdots x_k^{dn/k+a_k}$
in the generating function
$\exp(\sum_{1 \leq j < l \leq k} x_j x_l)$. Then as $n \to \infty$ we have $C_{a_1, a_2,\ldots, a_k}\sim C(s)$, where
\[
C(s) = (2\pi)^{-k} \left(\frac{k(k-1)}{dn}\right)^{(dn + s)/2}
2 e^{dn/2} (2\pi)^{k/2} \left(\frac{k(k-1)}{dn}\right)^{k/2}
(2k-2)^{-1/2} (k-2)^{-(k-1)/2}.
\]
\end{lemma}
Hence we deduce
\bel{eyLong}
\left| \Pnd \right| \ex (Y) \sim
\binom{n}{n/k, n/k, \ldots, n/k} (dn/k)!^k C(0).
\ee
Combining this with the well-known formula for the number of pairs on $dn$ points,
\bel{eq:Pnd}
|\Pnd| = (dn-1)!! = \frac{(dn)!}{(dn/2)!2^{dn/2}}\sim \sqrt{2}\left(\frac{dn}{e}\right)^{dn/2},
\ee
and after some basic manipulations using Stirling's 
formula we obtain the estimate for $\ex (Y)$ stated in the proposition.

Next we estimate the expected value of $YX_m$ where $Y$ is the number of
balanced $k$-colourings and $X_m$ the number of length-$m$ cycles.
It is more convenient to count rooted oriented cycles, which 
introduces a factor
of $2m$ into our calculations.  It will be helpful to have
the following definitions. For a rooted oriented cycle in a coloured
graph, define its
\emph{colour type} to be the sequence $T$ of colours on its vertices.
For $j = 1, 2, \ldots, k$, let $\alpha_j(T)$ denote the number of 
vertices in $T$ which have colour $j$. Note that the sum $\sum_j 
\alpha_j(T)$ is $m$.

To calculate the
expected value of $YX_m$, we will count, for each balanced
$k$-colouring and each
rooted oriented $m$-cycle, the number of pairings which contain this cycle
and respect this colouring.

As before, there are $\binom{n}{n/k, n/k, \ldots, n/k}$ ways to choose
the balanced $k$-colouring. All are equivalent, so fix one. To enumerate the
cycles and pairings which respect this colouring, we will sum over all colour
types $T$. Once a colour type has been chosen, each vertex of the
cycle can be placed
in the pairing model by choosing a vertex of the correct colour and
an ordered pair of
points in that vertex to be used by the cycle. Hence, in total, there
are asymptotically
$\left( d(d-1)n/k\right)^m$ ways to place the rooted oriented cycle
in the pairing model. We now have
$$
\ex (YX_m) \sim
\frac{1}{2m}
\binom{n}{n/k, n/k, \ldots, n/k}
\left(\frac{ d(d-1)n}{k}\right)^m
\frac{1}{\left| \Pnd \right|}
\sum_T f(T),
$$
where $f(T)$ is the number of pairings which respect a fixed balanced 
$k$-colouring and
fixed rooted oriented cycle of colour type $T$.
To count these pairings, suppose there are $b_{ij}=b_{ji}$ edges between
colour class $i$ and colour class $j$ (for $1\le i, j\le k$ and $i\ne j$), excluding the edges
of the prescribed cycle. The colours of the
neighbours of all of the unmatched points of colour class $i$
can be then chosen in
$(dn/k - 2\alpha_i(T))!/\prod_{1\le j\le k\atop i\ne j} b_{ij}!$ ways. 
After this determination
is made, edges are constructed by putting a perfect matching between
the corresponding points in each pair of classes, in one of 
$\prod_{1\le i<j\le k}b_{ij}!$
ways. Thus we have
\begin{eqnarray*}
f(T) &=&
\sum_{\{b_{ij}\}}
\left( \prod_{i=1}^k \frac{(dn/k - 2\alpha_i(T))!}{\prod_{1\le j\le k\atop i\ne j} b_{ij}!} \right)
\prod_{1\le i<j\le k}b_{ij}! \\
&=&
\sum_{\{b_{ij}\}}
\frac{\prod_i(dn/k - 2\alpha_i(T))!}{\prod_{i<j}b_{ij}!} \\
&\sim&
\frac{(dn/k)!^k}{(dn/k)^{2m}}
\sum_{\{b_{ij}\}}
\frac{1}{\prod_{1\le i<j\le k}b_{ij}!} \\
&\sim&
\frac{(dn/k)!^k}{(dn/k)^{2m}}
\Big[\prod_{l=1}^k x_l^{dn/k-2\alpha_l(T)}\Big] \prod_{1\le i<j\le k}\sum_{l\geq 0} 
\frac{(x_i x_j)^l}{l!} \\
&\sim&
\frac{(dn/k)!^k}{(dn/k)^{2m}}
\Big[\prod_{l=1}^k x_l^{dn/k-2\alpha_l(T)}\Big] \exp\left(\sum_{1\le i<j\le k} x_i x_j\right).
\end{eqnarray*}
By Lemma~\ref{lem:C(s)} the asymptotic value of the coefficient in the
last expression is $C(-2m)$, making the entire expression independent 
of $T$.
Moreover, the number $t_m$ of possible colour types for a rooted oriented cycle of length $m$ satisfies the obvious recurrence $t_m + t_{m-1} = k(k-1)^{m-1}$ with $t_1=0$. So we have $t_m=(k-1)^m + (k-1)(-1)^m$ and therefore
\begin{eqnarray*}
\ex (YX_m) \sim
\frac{1}{2m}
\binom{n}{n/k, n/k, \ldots, n/k}
\left(\frac{ d(d-1)n}{k}\right)^m
\frac{1}{\left| \Pnd \right|}
\frac{(dn/k)!^k}{(dn/k)^{2m}}
\\
\times
((k-1)^m + (k-1)(-1)^m)
C(-2m).
\end{eqnarray*}
Comparing this expression with $\eqn{eyLong}$ we see that
\begin{eqnarray*}
\frac{ \ex (YX_m) }{ \ex (Y) }
&\sim&
\frac{1}{2m}
\left(\frac{ d(d-1)n}{k}\right)^m
((k-1)^m + (k-1)(-1)^m)
\frac{1}{(dn/k)^{2m}}
C(-2m)/C(0)
\\
&\sim&
\frac{1}{2m}
\left(\frac{ d(d-1)n}{k}\right)^m
((k-1)^m + (k-1)(-1)^m)
\frac{1}{(dn/k)^{2m}}
\left(\frac{ k(k-1) }{ dn } \right)^{-m}
\\
&\sim&
\frac{(d-1)^m}{2m}
\left(1 + \frac{(-1)^m}{(k-1)^{m-1}}\right).
\\
&\sim&
\lambda_m
(1 + \delta_m).
\end{eqnarray*}
The above argument is easily extended to work for higher moments, by counting the pairings that contain a given balanced
$k$-colouring and set of oriented cycles of the appropriate lengths. The contribution from cases where the cycles intersect turns out to be negligible, for the following reasons. Suppose that the cycles form a subgraph $H$ with $\nu$ vertices and $\mu$ edges, and the total length of cycles is $\nu_0$. Then in the case of disjoint cycles, $\nu=\mu=\nu_0$. A factor of $\Theta(n^{\nu_0-\nu})$ is lost if there is a reduction in the number of vertices of $H$, compared with the disjoint case, because of the reduced number of ways of placing the cycles on the coloured vertices. Similarly, a factor $\Theta(n^{\nu_0-\mu})$ is gained in the function $f$ for the reduction in the number of edges of $H$, because of the corresponding increase in the number of points to be paired up at the end.   
Thus, the contribution from such an arrangement of cycles to the quantity being estimated is of the order of $n^{\nu-\mu}$ times that of the contribution from disjoint cycles.  In all non-disjoint cases, $H$ has more edges than vertices, since its minimum degree is at least $2$, and it has at least one vertex of degree at least $3$. There are only finitely many isomorphism types of $H$ to consider, so the contribution from the case of disjoint cycles is of the order of $n$ times the rest. The significant terms in this case decompose into a product of the factors corresponding to the individual cycles, and we obtain
$$
\ex(Y[X_1]_{p_1}\cdots [X_j]_{p_j}) / \ex (Y) \sim
\prod_{m=1}^j \left( \lambda_m \left( 1 + \delta_m \right)\right)^{p_m}
$$
as claimed.
\qed
It only remains to prove Lemma~\ref{lem:C(s)}. Before doing so, we need the following result, which will be used several times in the paper.
\begin{lemma}
\lab{lGauss}
Let $k$ be a positive integer. Define the function $f:\real^k\to\BC$ by
\[
f(\vtheta) = ia(n, \vtheta) -c_1 n \vtheta^\transpose B\vtheta
\]
where $i$ is the imaginary unit, 
$a$ is a real function, $B$ is a fixed $k$-by-$k$ positive
definite real matrix, and $c_1>0$ is a real constant. Let $\delta=c_2 n^{-1/2}\log n$ for some real constant $c_2>0$.
Then, as $n\to\infty$,
\[
\int_{[-\delta,\delta]^k} e^{f(\vtheta)} d\vtheta
=
\int_{[-\infty,\infty]^k} e^{f(\vtheta)} d\vtheta
+O(e^{-c(\log n)^2})
\]
for some constant $c>0$.
\end{lemma}
\proof
Note that in order to bound the integral of $e^{f(\theta)}$ over $[-\infty,\infty]^k\setminus [-\delta,\delta]^{k}$ we only need to consider the real part of $f(\theta)$. Moreover,
since $B$ is positive definite we have $\vtheta^\transpose B\vtheta \ge \lambda |\vtheta|^2$ where $\lambda>0$ is the smallest eigenvalue of $B$. The proof is elementary in view of these two observations.
\qed

\proofof{Lemma~\ref{lem:C(s)}}
We will use the saddlepoint method.
First we use Cauchy's formula to express $C_{a_1, a_2,\ldots, a_k}$ as an integral over
the product of circles $z_j = re^{i\theta_j}$, 
 $-\pi \leq \theta_j \leq \pi$
 for $j=1, 2, \ldots, k$,
where $r = \sqrt{dn/k(k-1)}$.
\begin{eqnarray*}
C_{a_1, a_2,\ldots, a_k} &=&
\frac{1}{(2\pi i)^k}
\int_{|z_1|=r}
\int_{|z_2|=r}
\cdots
\int_{|z_k|=r}
\frac{ \exp\left(\sum_{j<l} z_j z_l \right) }
{ z_1^{dn/k + a_1 + 1} z_2^{dn/k + a_2 + 1} \cdots
  z_k^{dn/k + a_k + 1}
}
dz_1 dz_2 \cdots dz_k
\\
&=&
\frac{1}{(2\pi)^k}
\int_{-\pi}^{\pi}
\int_{-\pi}^{\pi}
\cdots
\int_{-\pi}^{\pi}
\frac{ \exp\left(\sum_{j<l} (r e^{i\theta_j})(r e^{i\theta_l})\right) }
{ (re^{i\theta_1})^{dn/k + a_1 } (re^{i\theta_2})^{dn/k + a_2 } \cdots
  (re^{i\theta_k })^{dn/k + a_k }
}
  d\theta_1 d\theta_2 \cdots d\theta_k
\\
&=&
\frac{1}{(2\pi)^k r^{dn + s}}
\int_{-\pi}^{\pi}
\int_{-\pi}^{\pi}
\cdots
\int_{-\pi}^{\pi}
\frac{ \exp ( r^2 \sum_{j<l} e^{i(\theta_j + \theta_l)})}
      { \exp (i \sum_j (dn/k + a_j) \theta_j ) }
d\theta_1 \cdots d\theta_k.
\end{eqnarray*}
Let $g(\mathbf{\theta})$ denote the integrand in the last expression above.
Letting $\mathbf{1}$ denote the vector of 1's, consider the image of
$g(\mathbf{\theta})$ under the transformation
$\mathbf{\theta} \mapsto \mathbf{\theta} + \pi \mathbf{1}$. It is clear
that the numerator is fixed by this transformation. The denominator
becomes
\begin{eqnarray*}
\exp (i \sum_j (dn/k + a_j) (\theta_j  + \pi) )
&=&
\exp (i \sum_j (dn/k + a_j) \theta_j )
\exp (i (dn + s)  \pi )
\\
&=&
\exp (i \sum_j (dn/k + a_j) \theta_j )
\end{eqnarray*}
since $dn$ (the sum of the vertex degrees) and $s$ are both even. So
$g(\mathbf{\theta})$ is fixed by this transformation.
Letting $\delta = \log n / \sqrt{n}$, this means that the integrals of
$g(\mathbf{\theta})$ over regions
$\{ \mathbf{\theta} : |\theta_j| \leq \delta, j = 1, 2, \ldots, k \}$
and 
$\{ \mathbf{\theta} : \pi - \delta \leq |\theta_j| \leq \pi, j = 1, 2, \ldots, k \}$
are equal.
We will prove that the integral $I$ of $g(\mathbf{\theta})$ over each 
of these regions is asymptotically equal to
\begin{eqnarray*}
I &=&
e^{dn/2}
(2\pi)^{k/2}
\left(
\frac{k(k-1)}{dn}
\right)^{k/2}
(2k-2)^{-1/2}
(k-2)^{-(k-1)/2}
\\
&=&
K\exp\left( dn/2 -\frac{k}{2} \log n \right),
\end{eqnarray*}
where $K$ is a constant, and we will show that the integral over the 
remaining regions is asymptotically smaller.
 From these results the proposition follows.

To prove that the integral over vectors $\mathbf{\theta}$ in the 
remaining regions is asymptotically
smaller, there are two cases: either $|\theta_{j^*}| \leq \delta $ and
$\pi - \delta \leq |\theta_{l^*}| \leq \pi$ for some distinct $j^*$ and $l^*$, or
$\delta \leq |\theta_{j^*}| \leq \pi - \delta$ for some $j^*$.

In the first case, suppose that $|\theta_{j^*}| \leq \delta $ and
$\pi - \delta \leq |\theta_{l^*}| \leq \pi$ for some distinct $j^*$ and $l^*$. Then
$\pi - 2\delta \leq |\theta_{j^*} + \theta_{l^*}| \leq \pi + 2\delta$
and hence $\cos( \theta_{j*} + \theta_{l*} ) \leq 0$.
So
\begin{eqnarray*}
  |g(\mathbf{\theta})|
&=&
\exp \left( r^2 \sum_{j<l} \cos( \theta_j + \theta_l )  \right)
\\
&\leq&
\exp \left(  r^2 ( \binom{k}{2} - 1 ) + r^2 \cos( \theta_{j^*} + \theta_{l^*} )
\right)
\\
&\leq&
\exp \left( r^2 ( \binom{k}{2} - 1 )  \right)
\\
&=&
\exp \left( \frac{dn}{2} - \frac{dn}{k(k-1)} \right)
\\
&=&
o(I).
\end{eqnarray*}

In the second case we suppose that
$\delta \leq |\theta_{j^*}| \leq \pi - \delta$ for some $j^*$.
If there is a value of $l^*$ for which $|\theta_{j^*} + 
\theta_{l*}| > \delta / 2$ then
$\delta/2 < |\theta_{j^*} + \theta_{l*}| < 2\pi - \delta / 2$. This means
$$
\cos(\theta_{j^*} + \theta_{l*}) < \cos( \delta / 2 )
= 1 - \frac{\delta^2}{8} + O(\delta^4)
$$
and hence
\begin{eqnarray*}
|g(\mathbf{\theta})|
&=&
\exp\left( r^2 \sum_{j<l} \cos( \theta_j + \theta_l ) \right)
\\
&\leq&
  \exp\left( r^2 ( \binom{k}{2} - 1 ) + r^2 \cos( \theta_{j^*} + 
\theta_{l^*} ) \right)
\\
&=&
  \exp\left( r^2 ( \binom{k}{2} - 1 ) + r^2 ( 1 - \frac{\delta^2}{8} + 
O(\delta^4) ) \right)
\\
&=&
  \exp\left( r^2  \binom{k}{2}  - r^2\frac{\delta^2}{8} + 
O(r^2\delta^4)  \right)
\\
&=&
  \exp\left( \frac{dn}{2} - \frac{d(\log n)^2}{8k(k-1)} + o(1) \right)
\\
&=&
o(I).
\end{eqnarray*}

Otherwise, there is no such $l^*$. That is, for all $l^*$ not equal to
$j^*$ we have $|\theta_{l^*} - (-\theta_{j^*})| \leq \delta/2$. This implies
that all $\theta_l$ with $l \not = j^*$ have the same sign and satisfy
$\delta/2  \leq |\theta_{l}| \leq \pi - \delta/2$. Since $k\geq 3$ we can
choose two distinct such $l$, say $l^*$ and $l^{**}$, and deduce
$$
\delta \leq |\theta_{l^*} + \theta_{l^{**}}| \leq 2\pi - \delta.
$$
Using the same argument as above, it follows that
$|g(\mathbf{\theta})| =  o(I)$.

This completes the proof that the integral of $g(\mathbf{\theta})$ 
over these regions is asymptotically negligible, as claimed.

It remains to show the integral of $g(\mathbf{\theta})$ over the region
$\{ \mathbf{\theta} : |\theta_j| \leq \delta, j = 1, 2, \ldots, k \}$ is
asymptotically equal to $I$. We begin by expanding
\begin{eqnarray*}
\log g(\mathbf{\theta})
&=&
r^2\binom{k}{2}
+ i(r^2(k-1)-dn/k + O(1)) \sum_{j=1}^k \theta_j
-\frac{1}{2}r^2 \sum_{j<l}(\theta_j + \theta_l)^2
+ O(r^2 \sum_{j=1}^k |\theta_j|^3)
\\
&=&
r^2\binom{k}{2}
-\frac{1}{2}r^2 \sum_{j<l}(\theta_j + \theta_l)^2
+o(1)
\end{eqnarray*}
since $r^2=dn/(k(k-1))$ and $|\theta_j| \leq \delta = \log n / \sqrt{n}$ for
all $j$.
The quadratic order term can be written as
$-\frac{1}{2}r^2 \sum_{j<l}(\theta_j + \theta_l)^2 = -\frac{1}{2} \mathbf{\theta^\transpose}A\mathbf{\theta}$. Here, $\theta^\transpose$ denotes the transpose of the column vector $\theta$ and $A$
is the matrix $A = r^2(\mathbf{1}\mathbf{1^\transpose}+(k-2)I_k)$, where $I_k$ is 
the $k$-by-$k$ identity matrix. By Lemma~\ref{lGauss} and since $A$ is positive definite, we have for some constant $c>0$
\begin{eqnarray*}
&&
\int_{-\delta}^{\delta}
\int_{-\delta}^{\delta}
\cdots
\int_{-\delta}^{\delta}
\exp\left(-\frac{1}{2}r^2 \sum_{j<l}(\theta_j + \theta_l)^2\right)
d\theta_1 d\theta_2 \cdots d\theta_k
\\
&=&
\int_{-\infty}^{\infty}
\int_{-\infty}^{\infty}
\cdots
\int_{-\infty}^{\infty}
\exp\left(-\frac{1}{2} \mathbf{\theta^\transpose}A\mathbf{\theta}\right)
d\theta_1 d\theta_2 \cdots d\theta_k
+ O(e^{-c\log^2 n}).
\end{eqnarray*}
It is well-known (see Equation 4.6.3 
in~\cite{DeB}) that
such integrals have the value $(2\pi)^{k/2} ({\rm det}A)^{-1/2}$, giving us
\begin{eqnarray*}
&&
\int_{-\infty}^{\infty}
\int_{-\infty}^{\infty}
\cdots
\int_{-\infty}^{\infty}
\exp\left(-\frac{1}{2}r^2 \sum_{j<l}(\theta_j + \theta_l)^2\right)
d\theta_1 d\theta_2 \cdots d\theta_k
\\
&=&
(2\pi)^{k/2} \left( r^{2k} (2k-2)(k-2)^{k-1}\right)^{-1/2}.
\end{eqnarray*}
We conclude
\begin{eqnarray*}
\int_{-\delta}^{\delta}
\int_{-\delta}^{\delta}
\cdots
\int_{-\delta}^{\delta}
g(\mathbf{\theta})
d\theta_1 d\theta_2 \cdots d\theta_k
&\sim&
e^{r^2\binom{k}{2}}
(2\pi)^{k/2} \left( r^{2k} (2k-2)(k-2)^{k-1}\right)^{-1/2}.
\\
&\sim&
e^{dn/2}
(2\pi)^{k/2} \left( \frac{k(k-1)}{dn} \right)^{k/2} 
(2k-2)^{-1/2}(k-2)^{-(k-1)/2}
\\
&=&
I
\end{eqnarray*}
as claimed.
\qed


\section{Second moment: proof of Proposition~\ref{prop:moments}$\mbf{(b)}$}
Throughout this section, fix positive integers $d$ and $k\ge 3$ satisfying $d<2(k-1)\log(k-1)$. Assume $2$ divides $dn$ and $k$ divides $n$.
Let $C_1$ and $C_2$ be balanced $k$-colourings of a pairing $P\in\Pnd$. 
The \emph{colour count} of $(C_1,C_2)$ is the $k$-by-$k$ matrix $M=[m_{p,q}]$
where $m_{p,q}n/k$ is the number of cells coloured $p$ in $C_1$ and coloured
$q$ in $C_2$. 
Let $\M$ be the
set of $k$-by-$k$ \emph{doubly stochastic} matrices (i.e.~nonnegative matrices with each row sum and column sum 
equal to $1$).
Since the colourings are balanced, we must have $M\in\M$.
Define $T(M)$ to be the set of 
triples $(P,C_1,C_2)$ where
$P\in\Pnd$ and $(C_1, C_2)$ is a pair of balanced $k$-colourings 
of $P$ having colour count $M$.
  Then,
\bel{eSum}
\ex(Y^2) = \sum_{M\in\M\cap\frac{k}{n}\Z^{k^2}}\frac{|T(M)|}{|\Pnd|}.
\ee
In order to estimate the sum in~\eqref{eSum}, we first obtain an exact expression for $|T(M)|$, where $M=[m_{p,q}]$ is any $k$-by-$k$ doubly stochastic matrix whose entries are integer multiples of $k/n$.

For all $1\le p,q \le k$ we must choose $m_{p,q}n/k$ cells to be assigned the
colour $p$ in the first colouring and $q$ in the second colouring. We say that
such a cell and its points have label $(p,q)$.
The number of ways of doing this is given by

the multinomial coefficient
\[
\frac{n!}{\prod\limits_{1\le p,q\le k} (m_{p,q}n/k)!}.
\]
Now we must select the edges of the pairing in a way which is compatible with the two colourings.
Suppose we know the number $b_{pqrs}$ of edges from points 
labelled $(p,q)$ to points labelled $(r,s)$ for all $1\le p,q,r,s\le k$
with $p\ne r$ and $q\ne s$. Then we choose,
 for each ordered pair of labels $((p,q),(r,s))$, which
$b_{pqrs}$ of the points labelled $(p,q)$ will be paired with points 
labelled $(r,s)$. The number of ways of doing this is
\[
\prod_{1\le p,q\le k}\frac{(dm_{p,q}n/k)!}
{\prod\limits_{1\le r,s\le k\atop r\ne p,s\ne q} b_{pqrs}!}.
\]
Finally, for each unordered pair 
of labels $\{(p,q),(r,s)\}$, we choose a bijection between the points 
labelled $(p,q)$ and the points 
labelled $(r,s)$ that were designated to be paired with each other.
The number of ways of doing this is
\[
\prod_{1\le p,q,r,s\le k\atop p<r,q\ne s} b_{pqrs}!.
\]
Observe that the only restrictions on $b_{pqrs}$ required in our counting are
\begin{align}
b_{pqrs}&=b_{rspq} & \forall p,q,r,s\in\{1,\ldots,k\}, r\ne p, s\ne q,
\label{eq:B1}
\\
\sum_{1\le r,s\le k\atop r\ne p,s\ne q} b_{pqrs} &= d m_{p,q} n/k & \forall p,q\in\{1,\ldots,k\}.
\label{eq:B2}
\end{align}
Thus, the total number of triples $(P,C_1,C_2)$ in $T(M)$ is
\bel{eTM1}
|T(M)|=n!
\left(\prod_{1\le p,q\le k} \frac{(dm_{p,q}n/k)!}{(m_{p,q}n/k)!}\right)
\sum_{\cB(M,n)} \prod_{1\le p,q,r,s\le k\atop p<r,q\ne s} \frac{1}{b_{pqrs}!},
\ee
where $\cB(M,n)$ is the set of tuples of non-negative integers $(b_{pqrs})_{1\le p,q,r,s\le k\atop r\ne p, s\ne q}$ satisfying~\eqref{eq:B1} and~\eqref{eq:B2}. Note that~\eqref{eTM1} is the expression used in~\cite{AM}.

As we shall see later, the main weight of the sum in~\eqref{eSum} corresponds to terms in which $M$ is `near' $(1/k)J_k$, where $J_k$ denotes the $k$-by-$k$ matrix of ones.
To state this more precisely, for $\delta>0$ and any positive integer $p$, we define $\B_p(\delta)$ to be the set of $p$-by-$p$ matrices $M=[m_{ij}]$ for which $\max_{i,j}|m_{ij}-(1/k)|<\delta$.
It will be shown that the essential contribution to~\eqref{eSum} comes from terms such that $M\in \B_k(\epsilon\log n/n^{1/2})$, where $\epsilon>0$ is a small constant that will be specified later.
Thus we now proceed to bound $|T(M)|/|\Pnd|$ for each $M\in\M\cap(k/n)\Z^{k^2}$, and then find more accurate asymptotic expressions for $M\in \B_k(\epsilon\log n/n^{1/2})$.

We begin by estimating the ratio of factorials in~\eqref{eTM1}.
 Recall that one can write Stirling's formula as
$x!=\xi(x)(x/e)^x$ where $\xi$ is a function satisfying 
$\xi(x)\sim \sqrt{2\pi x}$ as $x\to\infty$ and
$\xi(x)\ge1$ for all $x\ge 0$. Thus,
\begin{eqnarray}
\prod_{1\le p,q\le k}
\frac{(dm_{p,q}n/k)!}{(m_{p,q}n/k)!}
&=&
\prod_{1\le p,q\le k}
\frac{\xi(dm_{p,q}n/k) (dm_{p,q}n/(ke))^{dm_{p,q}n/k}}
     {\xi( m_{p,q}n/k) ( m_{p,q}n/(ke))^{ m_{p,q}n/k}}
\notag\\
&=&
\prod_{1\le p,q\le k}
\frac{\xi(dm_{p,q}n/k)}{\xi( m_{p,q}n/k)}
d^{dm_{p,q}n/k}
\left(\frac{m_{p,q}n}{ke}\right)^{(d-1)m_{p,q}n/k}
\notag\\
&=&
d^{dn}
\left(\frac{n}{ek}\right)^{(d-1)n}
\prod_{1\le p,q\le k}
\frac{\xi(dm_{p,q}n/k)}{\xi( m_{p,q}n/k)}
{m_{p,q}}^{(d-1)m_{p,q}n/k},
\label{eq:FacRat1}
\end{eqnarray}
where in the final step we used $\sum_{1\le p,q\le k} m_{p,q}=k$ 
which holds because $M$ is doubly stochastic. (Throughout the article we use the convention $0^0=1$ and $0\log0=0$.)
Moreover, since $\xi(m_{p,q}n/k)\ge1$ and $\xi(dm_{p,q}n/k)=O(n^{1/2})$ for each $m_{p,q}$, we obtain the following bound, which does not depend on the particular $M$.
\bel{eq:FacRat2}
\prod_{1\le p,q\le k}
\frac{(dm_{p,q}n/k)!}{(m_{p,q}n/k)!} =
O\left(n^{k^2/2}\right) d^{dn} \left(\frac{n}{ek}\right)^{(d-1)n}
\left(\prod_{1\le p,q\le k}
{m_{p,q}}^{m_{p,q}}\right)^{(d-1)n/k}.
\ee
Next we bound the inner sum in~\eqref{eTM1}, by further applying Stirling's formula, and also using the obvious crude bound $|\cB(M,n)|=O(n^{k^4})$ on the number of terms.
\begin{align*}
\sum_{\cB(M,n)}
\prod_{1\le p,q,r,s\le k\atop p<r,q\ne s} \frac{1}{b_{pqrs}!} 
&\sim
\sum_{\cB(M,n)}
\prod_{1\le p,q,r,s\le k\atop p<r,q\ne s} \frac{1}{\xi(b_{pqrs})(b_{pqrs}/e)^{b_{pqrs}}} 
\\
&\le
\sum_{\cB(M,n)}
\prod_{1\le p,q,r,s\le k\atop p<r,q\ne s} \frac{1}{(b_{pqrs}/e)^{b_{pqrs}}} 
\\
&=
O\left(n^{k^4}\right) e^{dn/2} \max_{\cB(M,n)}
\Bigg\{\prod_{1\le p,q,r,s\le k\atop p<r,q\ne s}\frac{1}{{b_{pqrs}}^{b_{pqrs}}}\Bigg\}.
\end{align*}
The following result allows us to derive a more explicit bound, conveniently expressed in terms of $M$ and $n$. The proof follows the ideas in~\cite{AM} and is given later in this section.
\begin{lemma}\label{lem:regtognm}
Let $M=[m_{pq}]$ be a fixed matrix in $\cM$, and let $(b_{pqrs})_{1\le p,q,r,s\le k\atop p<r,q\ne s}$ be any tuple of non-negative reals satisfying~\eqref{eq:B1} and~\eqref{eq:B2}. Then,
\[
\prod_{1\le p,q,r,s\le k\atop p<r,q\ne s}
\frac{1}{{b_{pqrs}}^{b_{pqrs}}} \le
\left(\frac{1}{\prod_{1\le p,q\le k} {m_{p,q}}^{m_{p,q}}}\right)^{dn/k} \left(\frac{\sum_{1\le p,q,r,s\le k\atop p\ne r,q\ne s} m_{p,q}m_{r,s}}{dn}\right)^{dn/2}.
\]
\end{lemma}
Hence, we immediately deduce that
\bel{eq:InnerSum1}
\sum_{\cB(M,n)}
\prod_{1\le p,q,r,s\le k\atop p<r,q\ne s} \frac{1}{b_{pqrs}!} 
=
O\left(n^{k^4}\right) \left(\frac{1}{\prod_{1\le p,q\le k} {m_{p,q}}^{m_{p,q}}}\right)^{dn/k} \left(\frac{\sum_{1\le p,q,r,s\le k\atop p\ne r,q\ne s} m_{p,q}m_{r,s}}{dn/e}\right)^{dn/2}.
\ee
Finally, define for any $M\in\cM$
\begin{align}
\varphi(M) &= -\frac{1}{k}\sum_{1\le p,q\le k} m_{p,q}\log m_{p,q} + \frac{d}{2}\log\left(\frac{1}{k^2}\sum_{1\le p,q,r,s\le k\atop p\ne r,q\ne s} m_{p,q}m_{r,s}\right)
\notag\\
&= -\frac{1}{k}\sum_{1\le p,q\le k} m_{p,q}\log m_{p,q} + \frac{d}{2}\log\left(1-\frac{2}{k}+\frac{1}{k^2}\sum_{1\le p,q\le k} {m_{p,q}}^2\right).
\label{eq:G}
\end{align}
(Recall the convention $0^0=1$ and $0\log0=0$.)
By combining~\eqref{eq:FacRat2}, \eqref{eq:InnerSum1}, the Stirling-formula estimate  
$n!\sim\sqrt{2\pi n}(n/e)^n$ and~\eqref{eq:Pnd}, we have
\begin{align}
\frac{|T(M)| }{|\Pnd|} &= \frac{\sqrt{2\pi n}(n/e)^n}{\sqrt{2}(dn/e)^{dn/2}}\,
O\left(n^{k^2/2}\right) d^{dn} \left(\frac{n}{ek}\right)^{(d-1)n}
\left(\prod_{1\le p,q\le k} {m_{p,q}}^{m_{p,q}}\right)^{(d-1)n/k}
\notag\\
&\quad \times O\left(n^{k^4}\right) \left(\frac{1}{\prod_{1\le p,q\le k} {m_{p,q}}^{m_{p,q}}}\right)^{dn/k} \left(\frac{\sum_{1\le p,q,r,s\le k\atop p\ne r,q\ne s} m_{p,q}m_{r,s}}{dn/e}\right)^{dn/2}
\notag\\
&= O\left(n^{k^4+k^2/2+1/2}\right) k^n
\left(\prod_{1\le p,q\le k}
\frac{1}{{m_{p,q}}^{m_{p,q}}}\right)^{n/k}
\left(\frac{1}{k^2}\sum_{1\le p,q,r,s\le k\atop p\ne r,q\ne s} m_{p,q}m_{r,s}\right)^{dn/2}
\notag\\
&\le \poly(n) k^n e^{n\varphi(M)},
\label{eq:term}
\end{align}
for some polynomial $\poly(n)$ not depending on the particular $M$.
This is a sufficient bound on $|T(M)|/|\Pnd|$, since it will allow us to show that the sum in~\eqref{eSum} receives a negligible contribution due to terms with $M$ away from $(1/k)J_k$.

It remains to find a suitable asymptotic expression of $|T(M)|/|\Pnd|$ for $M\in\B_k(\epsilon\log n/n^{1/2})$, by improving our previous estimates in~\eqref{eq:FacRat2} and~\eqref{eq:InnerSum1}.
Before that, we state two technical algebraic results which will be needed in the asymptotic calculations.
Hereinafter, $\onesvec$ represents the $k$-dimensional vector of ones,
while $\onesvec^{(i)}$ represents the $k$-dimensional vector with entry $1$ at position $i$ and $0$ elsewhere; $I_k$ denotes the $k$-by-$k$ identity matrix; $\mvec A$ is the vector formed by stacking the columns of a matrix $A$ to form a single column vector; and $A^{\otimes 2}$ is simply $A\otimes A$, with the standard notation $\otimes$ for the Kronecker product of matrices.
\begin{lemma}
\lab{lem:Evec2}
Consider the vectors
\[
\begin{array}{ll}
f^{(p)} = \frac{\sqrt{p}}{\sqrt{p+1}}\left(\frac{-1}{p}\sum_{l=1}^p\onesvec^{(l)}+\onesvec^{(p+1)}\right), & 1\le p \le k-1 \\
f^{(k)} = \frac{1}{\sqrt{k}}\onesvec. &
\end{array}
\]
Define $f^{(p,q)}=f^{(p)}\otimes f^{(q)}$ for $1\le p,q \le k$.
\begin{description}
\item{(a)}
An orthonormal basis of eigenvectors for the matrix 
$(J_k-I_k)^{\otimes 2}+(k-1)^2I_{k^2}$ is given by $\{f^{(p,q)}\}_{p,q=1}^k$ with corresponding eigenvalues
\[
\begin{array}{ll}
\lambda_{p,q} = k^2-2k+2,& 1\le p,q \le k-1 \\
\lambda_{p,k} = (k-1)(k-2),& 1\le p \le k-1 \\
\lambda_{k,q} = (k-1)(k-2),& 1\le q \le k-1 \\
\lambda_{k,k} = 2(k-1)^2.&
\end{array}
\]
The smallest of these eigenvalues is $(k-1)(k-2)$.
\item{(b)}
Similarly, the eigenvectors of $(J_k+I_k)^{\otimes 2}$ are also $\{f^{(p,q)}\}_{p,q=1}^k$, and the corresponding eigenvalues are
1 with multiplicity $(k-1)^2$,
$k+1$ with multiplicity $2(k-1)$,
and $(k+1)^2$ with multiplicity 1.
\end{description}
\end{lemma}
\proof Immediate by checking that the eigenvectors satisfy the required properties.
\qed
\begin{lemma}
\lab{lem:Evec3}
Let $A=[a_{i,j}]$ be a $k$-by-$k$ matrix whose rows and columns each
have sum $0$. Define $\tilde A$ to be the submatrix formed 
from $A$ by deleting the last row and column. 
Let $\{f^{(i,j)}\}_{i,j=1}^k$ be the orthonormal basis defined in the statement
of Lemma~\ref{lem:Evec2}. Then,
\begin{description}
\item{(a)}
$(\mvec A)^\transpose f^{(i,k)}=(\mvec A)^\transpose f^{(k,j)}=0$
for $1\le i,j \le k$,
\item{(b)}
$\displaystyle
\sum_{i=1}^{k-1}\sum_{j=1}^{k-1}\left((\mvec A)^\transpose f^{(i,j)}\right)^2
=
\sum_{i=1}^k \sum_{j=1}^k a_{i,j}^2,$  and
\item{(c)}
$\displaystyle
(\mvec \tilde A)^\transpose (J_{k-1}+I_{k-1})^{\otimes 2} \mvec \tilde A
=
\sum_{i=1}^k \sum_{j=1}^k a_{i,j}^2.
$
\end{description}
\end{lemma}
The proof of Lemma~\ref{lem:Evec3} is given later in this section.
We now proceed with the asymptotic calculations for $M\in \cM\cap(k/n)\Z^{k^2}\cap\B_k(\epsilon\log n/n^{1/2})$. So define the matrix $A=A(M)=[a_{p,q}]$ by $A=M-(1/k)J_k$. Note that each row and column of $A$ must
have sum $0$, and moreover $a_{p,q}<\epsilon\log n/n^{1/2}$ for each $p,q\in\{1,\ldots,k\}$. Let $\tilde A$ be the submatrix formed from $A$ by deleting the last row and column. Under the new assumptions, we first derive a new asymptotic formula for the expression computed in~\eqref{eq:FacRat1}.
For $m_{p,q}=1/k+a_{p,q}$ with $a_{p,q}=O(\log n/n^{1/2})$ we have
\[
\frac{\xi(dm_{p,q}n/k)}{\xi( m_{p,q}n/k)}
\sim
\frac{\sqrt{dm_{p,q}n/k}}{\sqrt{m_{p,q}n/k}}
\sim \sqrt{d}
\]
for $1\le p,q \le k$, and we expand
\begin{eqnarray*}
\sum_{1\le p,q\le k} m_{p,q}\log m_{p,q}
&=&
\sum_{1\le p,q\le k} \left(\frac{1}{k}+a_{p,q}\right)\left(\log\frac{1}{k}+\log\left(1+ka_{p,q}\right)\right)\\
&=&
\sum_{1\le p,q\le k} \left(-\frac{1}{k}\log k+\frac{k}{2}a_{p,q}^2+O(a_{p,q}^3)\right)\\
&=&
k\left(-\log k+\frac{1}{2}
(\mvec \tilde A)^\transpose (J_{k-1}+I_{k-1})^{\otimes 2} \mvec \tilde A
\right)+O(\log^3 n/n^{3/2})
\end{eqnarray*}
where we used Lemma~\ref{lem:Evec3}(c) to rewrite $\sum_{p,q}a_{p,q}^2$
in the final step.
Combining these estimates we can rewrite~\eqref{eq:FacRat1} as
\bel{eFactRat}
\prod_{1\le p,q\le k}
\frac{(dm_{p,q}n/k)!}{(m_{p,q}n/k)!}
\sim
d^{k^2/2}
d^{dn} \left(\frac{n}{ek^2}\right)^{(d-1)n}
\exp\left(n\frac{d-1}{2}
(\mvec \tilde A)^\transpose (J_{k-1}+I_{k-1})^{\otimes 2} \mvec \tilde A\right).
\ee
Next we rewrite the inner sum in~\eqref{eTM1} in terms of the natural generating function, letting square brackets denote the extraction of a coefficient
\begin{eqnarray*}
\sum_{\cB(M,n)}
\prod_{1\le p,q,r,s\le k\atop p<r,q\ne s} \frac{1}{b_{pqrs}!}
&=&
\Big[\prod_{1\le p,q\le k} {x_{p,q}}^{dm_{p,q}n/k}\Big]
\prod_{1\le p,q,r,s\le k\atop p<r,q\ne s}
\sum_{i=0}^\infty \frac{(x_{p,q}x_{r,s})^i}{i!}\\
&=&
\Big[\prod_{1\le p,q\le k} {x_{p,q}}^{dm_{p,q}n/k}\Big]
\exp\bigg(
 \frac{1}{2}
\sum_{1\le p,q,r,s\le k\atop p\ne r,q\ne s} 
x_{p,q}x_{r,s}\bigg).
\end{eqnarray*}
The following result provides an asymptotic characterisation of that coefficient. The proof uses the saddlepoint method, and can be found at the end of the section. Note that our still-unspecified $\epsilon$ is determined by the statement.
\begin{lemma}\label{lem:saddle}
There exists $\epsilon>0$ such that for each
$M\in \cM\cap(k/n)\Z^{k^2}\cap\B_k(\epsilon\log n/n^{1/2})$
the coefficient
$C$ of $\left[\prod_{1\le p,q\le k} {x_{p,q}}^{dm_{p,q}n/k}\right]$ in the generating function
$\exp\left(\frac{1}{2}\sum_{1\le p,q,r,s\le k\atop p\ne r,q\ne s} 
x_{p,q}x_{r,s}\right)$
satisfies
\[
C \sim \frac{\gamma(k)}{\sqrt\pi}
\left(\frac{1}{dn}\right)^{k^2/2}
\left(\frac{k(k-1)}{\sqrt{dn/e}}\right)^{dn}
\exp\left( 
\frac{-nd(k-1)^2 
(\mvec \tilde A)^\transpose (J_{k-1}+I_{k-1})^{\otimes 2} \mvec \tilde A
}{2(k^2-2k+2)}
\right),
\]
where $\gamma(k)$ is the constant
$\displaystyle
\gamma(k) = 
\frac{ k^{k^2} (k-1)^{k(k-1)} }
     { (2\pi)^{k^2/2-1/2} (k^2-2k+2)^{(k-1)^2/2} (k-2)^{k-1}}.$
\end{lemma}
Hence, in view of~\eqref{eFactRat} and Lemma~\ref{lem:saddle}, we obtain for $M\in \cM\cap(k/n)\Z^{k^2}\cap\B_k(\epsilon\log n/n^{1/2})$ an asymptotic expression for \eqref{eTM1} which improves the bound already stated in~\eqref{eq:term}
\bel{eq:asymp}
\frac{|T(M)| }{|\Pnd|} \sim \gamma(k)
\frac{(k-1)^{dn}}{n^{k^2/2-1/2} k^{(d-2)n}}
\exp \left(-n \frac{k^2-2k-d+2}{2(k^2-2k+2)}
(\mvec \tilde A)^\transpose (J_{k-1}+I_{k-1})^{\otimes 2} \mvec \tilde A
\right),
\ee
where $\gamma(k)$ is the constant defined in Lemma~\ref{lem:saddle}.

We are now in good shape to estimate the sum~\eqref{eSum}.
We begin by computing the contribution of the terms near $(1/k)J_k$.
For a $(k-1)$-by-$(k-1)$ 
matrix $M$ define $\overline M$ to be the $k$-by-$k$ matrix formed 
from $M$ by adding a new row and column so that every row sum 
and column sum is 1.
Recall the definition of $\B_p(\delta)$ and define $\overline{\B_p(\delta)}=\{\overline M \mid M\in \B_p(\delta)\}$.
Now we set
\[
\delta = \frac{\epsilon\log n}{(k-1)^2 n^{1/2}},
\]
and consider $M=\overline{M'}$ for $M'\in \B_{k-1}(\delta)\cap \frac{k}{n}\Z^{(k-1)^2}$.
A straightforward application of the triangle inequality shows that $\overline{\B_{k-1}(\delta)}\subseteq\B_k((k-1)^2\delta)$, and therefore
 $M \in \B_k(\epsilon\log n/n^{1/2})$. So $M$
is nonnegative (for large enough $n$) and hence $M\in \M$. Furthermore, the entries of 
$M$ are in $\frac{k}{n}\Z$ because they are integer linear
combinations of $k/n$ and $1=\frac{k}{n}\times \frac{n}{k}$, using the
fact that $k$ divides $n$. This shows that
\[
\frac{|T(M)|}{|\Pnd|}
\]
is a term in the sum \eqn{eSum}, suggesting that we express \eqn{eSum}
as $\ex(Y^2) = S_1+S_2$ where
\[
S_1 = \sum_{M'\in \B_{k-1}(\delta)\cap \frac{k}{n}\Z^{(k-1)^2}}
\frac{|T(\overline{M'})|}{|\Pnd|}
\]
and $S_2$ is the sum of the remaining terms.
Notice moreover that if $M'\in \B_{k-1}(\delta)\cap \frac{k}{n}\Z^{(k-1)^2}$ and $M=\overline{M'}$, then the matrix $A=M-(1/k)J_k$ has all entries $a_{p,q} < \epsilon n^{-1/2}\log n$. Hence the expansion given in~\eqref{eq:asymp} is valid, and we can express
\bel{eS1}
S_1 \sim \gamma(k) \frac{(k-1)^{dn}}{n^{k^2/2-1/2} k^{(d-2)n}} \sum_{M'\in \B_{k-1}(\delta)\cap \frac{k}{n}\Z^{(k-1)^2}} e^{nf(M')},
\ee
where
\[
f(M') = - \frac{k^2-2k-d+2}{2(k^2-2k+2)}
(\mvec \tilde A)^\transpose (J_{k-1}+I_{k-1})^{\otimes 2} \mvec \tilde A,
\]
and $\tilde A = M'-(1/k)J_{k-1}$.
By iterating the Euler-Maclaurin summation formula
(see~\cite{AS72}, p.~$806$), we have
\[
\sum_{M'\in \B_{k-1}(\delta)\cap \frac{k}{n}\Z^{(k-1)^2}} e^{nf(M')} \sim
\left(\frac{n}{k}\right)^{(k-1)^2}
\int_{M'\in \B_{k-1}(\delta)} e^{nf(M')}dM'.
\]
Setting $H = \frac{k^2-2k-d+2}{2(k^2-2k+2)} (J_{k-1}+I_{k-1})^{\otimes 2}$ and using Lemma~\ref{lem:Evec2} (b) for the eigenvalues of $(J_{k-1}+I_{k-1})^{\otimes 2}$, we deduce that $H$ is positive definite and has determinant
\begin{eqnarray}
k^{2k-2}\left(\frac{k^2-2k-d+2}{k^2-2k+2}\right)^{(k-1)^2}.
\lab{eDet}
\end{eqnarray}
Here we used the fact that $k^2-2k-d+2>0$, which is guaranteed from our assumptions on $d$ and $k$.
Now we may apply Lemma~\ref{lGauss} to conclude
\begin{align}
S_1 &\sim \gamma(k) \frac{(k-1)^{dn}}{n^{k^2/2-1/2} k^{(d-2)n}} \left(\frac{n}{k}\right)^{(k-1)^2}
\bigg(\int_{[-\infty,\infty]^{(k-1)^2}} e^{nf(M')}dM' + O\big(e^{-c\log^2n}\big) \bigg)
\notag\\
&\sim
\gamma(k) \frac{(k-1)^{dn}}{n^{k^2/2-1/2} k^{(d-2)n}}
\left(\frac{n}{k}\right)^{(k-1)^2}
\frac{(2\pi/n)^{(k-1)^2/2}}{|\mdet H|^{1/2}}
\notag\\
&=
\frac{k^{k} (k-1)^{k(k-1)}} {(k^2-2k-d+2)^{(k-1)^2/2} (2\pi(k-2))^{k-1}}
n^{-(k-1)} k^{2n} \left(1-\frac{1}{k}\right)^{dn}.
\label{eS1c}
\end{align}
To prove the proposition it suffices to show $S_2=o(S_1)$.
Let $M$ be an index of any term of $S_2$.
This implies $M\not\in \overline{ \B_{k-1}(\delta)}$, so we must have 
$M\in\M\setminus \B_k(\delta)$
since $\B_k(\delta)\cap\M\subseteq\overline{ \B_{k-1}(\delta)}$.
Recall now the definition of $\varphi$ in~\eqref{eq:G}, and obtain by direct substitution
\[
\varphi\left(\frac{1}{k}J_k\right) = \log k + d\log\left(1-\frac{1}{k}\right).
\]
{From} Theorem~7 in~\cite{AN} (see also (5) in the same paper) we have that if $d<d_{k-1}:=2(k-1)\log(k-1)$ for each $M\in\cM$,
\begin{align*}
\varphi(M) &\le \varphi\left(\frac{1}{k}J_k\right) - \frac{d_{k-1}-d}{4(k-1)^2}\left(\sum_{p,q} {m_{p,q}}^2 - 1\right)
\\
&= \log k + d\log\left(1-\frac{1}{k}\right) - \frac{d_{k-1}-d}{4(k-1)^2} \sum_{p,q} \left(m_{p,q}-\frac{1}{k}\right)^2.
\end{align*}
In particular, for each $M\in\M\setminus \B_k(\delta)$,
\begin{equation}
\varphi(M) \le \log k + d\log\left(1-\frac{1}{k}\right) - \frac{d_{k-1}-d}{4(k-1)^2} \delta^2.
\label{eq:boundG}
\end{equation}
By combining~\eqref{eq:boundG} with the bound on the general term obtained in~\eqref{eq:term} and also taking into account that the number of terms in $S_2$ is at most $O(n^{k^2})$, we conclude
\begin{align*}
S_2 &= O(n^{k^2}) \poly(n) k^{2n} \left(1-\frac{1}{k}\right)^{dn} \exp\left(- \frac{d_{k-1}-d}{4(k-1)^2} \epsilon\log^2 n\right)
\\&\le \poly'(n) k^{2n} \left(1-\frac{1}{k}\right)^{dn} n^{- \Theta(\log n)},
\end{align*}
for some polynomial $\poly'(n)$. Thus $S_2= o(S_1)$, and this completes the proof.
\qed
%
%
%
It only remains to prove Lemmas~\ref{lem:regtognm}, \ref{lem:Evec3} and~\ref{lem:saddle}.
\proofof{Lemma~\ref{lem:regtognm}}
Let $\cL=\cL(M)$ be the polytope consisting of all non-negative tuples $L=(\ell_{pqrs})_{1\le p,q,r,s\le k\atop p<r,q\ne s}$ in $\real^{k^2(k-1)^2/2}$ such that
\begin{equation}\label{eq:ell_m}
\sum_{1\le r,s\le k\atop r\ne p,s\ne q} \ell_{pqrs} = d m_{p,q}\qquad \forall p,q\in\{1,\ldots,k\},
\end{equation}
where for $p>r$ and $q\ne s$ we used the duplicate notation $\ell_{pqrs}=\ell_{rspq}$ to denote the coordinates.
For each $L\in\cL$, define
\[
\psi(L)=\prod_{1\le p,q\le k}
\frac{{m_{p,q}}^{dm_{p,q}}}{\prod_{1\le r,s\le k\atop p<r,q\ne s} {\ell_{pqrs}}^{\ell_{pqrs}}}.
\]
Our aim is to show that for all $L\in\cL$
\bel{eq:rescaled}
\psi(L) \le \left(\frac{\sum_{1\le p,q,r,s\le k\atop p\ne r,q\ne s} m_{p,q}m_{r,s}}{dk}\right)^{dk/2}.
\ee
This is indeed equivalent to the statement of the lemma, after setting $b_{pqrs}=\ell_{pqrs}n/k$ and performing straightforward manipulations.

As a first case, assume that $m_{pq}>0$ for all $p,q$ in $\{1,\ldots,k\}$.
We define $\hcL$ to be the polytope of all non-negative tuples $L=(\ell_{pqrs})_{1\le p,q,r,s\le k\atop p<r,q\ne s}$ in $\real^{k^2(k-1)^2/2}$ such that
\begin{equation}\label{eq:beta_sum}
\sum_{1\le p,q,r,s\le k\atop p<r,q\ne s} \ell_{pqrs} = dk/2,
\end{equation}
and for each $L\in\hcL$, let
\[
\widehat\psi(L) = \prod_{1\le p,q,r,s\le k\atop p<r,q\ne s} \left(\frac{m_{p,q}m_{r,s}}{{\ell_{pqrs}}}\right)^{\ell_{pqrs}}.
\]
Observe that $\cL\subset\hcL$, and moreover the restriction of $\widehat\psi$ to $\cL$ is equal to $\psi$.
Our goal is to maximise $\log\widehat\psi$ over $\hcL$, and thus provide a bound on $\psi$ over $\cL$. We first show that $\log\widehat\psi$ does not maximise on the boundary. Note that the boundary of $\hcL$ consists of points having some $0$ coordinate (but at least some coordinate must be strictly positive). Let us choose an arbitrary boundary point $L_0$, and assume without loss of generality that $\ell_{1122}=0$ and $\ell_{1221}>0$. Let $L_\epsilon$ be the point with the same coordinates as $L_0$ but replacing $\ell_{1122}$ by $\epsilon$ and $\ell_{1221}$ by $\ell_{1221}-\epsilon$. For small enough $\epsilon>0$, $L_\epsilon$ lies in $\hcL$, and moreover
\[
\lim_{\epsilon\to0^+}\frac{d}{d\epsilon}\log\widehat\psi(L_\epsilon)=+\infty,
\]
so the maximum cannot occur at $L_0$. Hence, $\log\widehat\psi$ is maximised in the interior of $\hcL$ and the maximiser(s) must satisfy the following Lagrange equations
\[
\log(m_{p,q}m_{r,s}) - \log\ell_{pqrs} - 1 = \lambda \qquad \forall 
p,q,r,s\in\{1,\ldots k\}, p<r, q\ne s.
\]
These are equivalent to
\[
\ell_{pqrs} = m_{p,q}m_{r,s} e^{-\lambda-1} \qquad \forall p,q,r,s\in\{1,\ldots k\}, p<r, q\ne s,
\]
which combined with~\eqref{eq:beta_sum} yield
\[
\ell_{pqrs} = dk \frac{m_{p,q}m_{r,s}}{\sum_{1\le p',q',r',s'\le k\atop p'\ne r',q'\ne s'} m_{p',q'}m_{r',s'}} \qquad \forall p,q,r,s\in\{1,\ldots k\}, p<r, q\ne s.
\]
{From} the uniqueness of the solution, we deduce that it must be the maximiser of $\log\widehat\psi$ (and $\widehat\psi$). The value of $\widehat\psi$ at this point can be easily computed by substitution
\[
\left(\frac{\sum_{1\le p,q,r,s\le k\atop p\ne r,q\ne s} m_{p,q}m_{r,s}}{dk}\right)^{dk/2},
\]
which proves the bound in~\eqref{eq:rescaled} under the assumption that $m_{pq}>0$ for all $p$ and $q$. To extend the argument to the other cases, we first define
\[
\cP = \{ (M,L)\st M\in\cM,\, L\in\cL(M) \},
\]
and with a mild abuse of notation denote by
\[
\psi(M,L) = \prod_{1\le p,q\le k} \frac{{m_{p,q}}^{dm_{p,q}}}{\prod_{1\le r,s\le k\atop p<r,q\ne s} {\ell_{pqrs}}^{\ell_{pqrs}}},
\]
the natural extension of $\psi$ to $\cP$. Notice that $\psi$ is continuous in $\cP$. We showed so far that, for any $M$ in the interior of $\cM$ (i.e.~$m_{pq}>0$) and any $L\in\cL(M)$,
\[
\psi(M,L) \le \left(\frac{\sum_{1\le p,q,r,s\le k\atop p\ne r,q\ne s} m_{p,q}m_{r,s}}{dk}\right)^{dk/2}.
\]
Hence, this inequality can be extended by continuity to any $M$ on the boundary of $\cM$, and thus to any $(M,L)\in\cP$. This concludes the proof of~\eqref{eq:rescaled}.
\qed


\proofof{Lemma~\ref{lem:Evec3}}
We begin by proving (a).
Let  $j\in\{1,2,\ldots,k\}$. 
Since $f^{(k,j)}=f^{(k)}\otimes f^{(j)}$ is a linear
combination of terms of the form 
$\onesvec\otimes\onesvec^{(q)}=\sum_{p=1}^k(\onesvec^{(p)}\otimes\onesvec^{(q)}), (1\le q \le k)$, we have that $(\mvec A)^\transpose f^{(k,j)}$ is a linear combination of
terms of the form 
$\sum_{p=1}^k(\mvec A)^\transpose (\onesvec^{(p)}\otimes\onesvec^{(q)})
=\sum_{p=1}^k a_{q,p}=0$ since the row sums of $A$ are $0$. 
 A similar argument shows $(\mvec A)^\transpose f^{(i,k)}=0$ for 
$i\in\{1,2,\ldots,k\}$
using the fact that the column sums of $A$ equal $0$.

To prove (b) we apply (a) to write
\[
\sum_{i=1}^{k-1}\sum_{j=1}^{k-1}\left((\mvec A)^\transpose f^{(i,j)}\right)^2
=
\sum_{i=1}^{k}\sum_{j=1}^{k}\left((\mvec A)^\transpose f^{(i,j)}\right)^2,
\]
which is the sum of the squares of the coordinates of 
$\mvec A$ in the basis given by $\{f^{(i,j)}\}_{i,j=1}^k$. Since the basis
is orthonormal, this expression is simply the square of the norm of $\mvec A$
with respect to the standard basis, 
$\sum_{i=1}^k \sum_{j=1}^k a_{i,j}^2$.

To prove part (c) we begin by writing
\begin{eqnarray*}
&&\sum_{i=1}^k\sum_{j=1}^k a_{i,j}^2\\
&=&
a_{k,k}^2
+\sum_{i=1}^{k-1}a_{i,k}^2
+\sum_{j=1}^{k-1}a_{k,j}^2
+\sum_{i=1}^{k-1}\sum_{j=1}^{k-1} a_{i,j}^2\\
&=&
\left( \sum_{i=1}^{k-1}\sum_{j=1}^{k-1}a_{i,j}  \right)^2
+\sum_{i=1}^{k-1}\left(-\sum_{j=1}^{k-1}a_{i,j}\right)^2
+\sum_{j=1}^{k-1}\left(-\sum_{i=1}^{k-1}a_{i,j}\right)^2
+\sum_{i=1}^{k-1}\sum_{j=1}^{k-1}a_{i,j}^2.
\end{eqnarray*}
Since
\[
\left( \sum_{i=1}^{k-1}\sum_{j=1}^{k-1}a_{i,j}  \right)^2
=
(\mvec \tilde A)^\transpose J_{k-1}^{\otimes 2}\mvec \tilde A,
\]
\[
\sum_{i=1}^{k-1}\left(\sum_{j=1}^{k-1}a_{i,j}\right)^2
= (\mvec \tilde A)^\transpose (J_{k-1}\otimes I_{k-1}) \mvec \tilde A,
\]
\[
\sum_{j=1}^{k-1}\left(\sum_{i=1}^{k-1}a_{i,j}\right)^2
= (\mvec \tilde A)^\transpose (I_{k-1}\otimes J_{k-1}) \mvec \tilde A,
\]
\[
\sum_{i=1}^{k-1}\sum_{j=1}^{k-1}a_{i,j}^2
= (\mvec \tilde A)^\transpose I_{k-1}^{\otimes 2}\mvec \tilde A,
\]
and $(J_{k-1}+I_{k-1})^{\otimes 2} = J_{k-1}^{\otimes 2}+(J_{k-1}\otimes I_{k-1})+(I_{k-1}\otimes J_{k-1})
+I_{k-1}^{\otimes 2}$, part (c) is proved.
\qed

Before proceeding with the proof of Lemma~\ref{lem:saddle}, we need the following technical result.
\begin{lemma}
\lab{lem:SaddlePair}
Let $\delta\in(0,2\pi/5)$ and fix an integer $k\ge 3$. For each
$1\le p,q \le k$, let
$-\pi \le \theta_{p,q} \le \pi$. Suppose 
$\max_{p,q} |\theta_{p,q}| > \delta$ and
$\min_{p,q} |\theta_{p,q}| < \pi - \delta$. 
Then there exist $p$, $q$, $r$, and $s$
with $p\ne r$ and $q\ne s$ such that 
\[
\frac{\delta}{2} \le |\theta_{p,q}+\theta_{r,s}| \le 2\pi - \frac{\delta}{2}.
\]
\end{lemma}
\proof
There are two cases. In the first case, suppose 
$\delta < |\theta_{p,q}| < \pi  - \delta$ for some $p$ and $q$.
Let $S$ be the set of pairs
\[
S= \{ (r,s) \mid 1\le r\le k, r\ne p, 1\le s\le k, s\ne q\}.
\]
The set $S$ is nonempty as $k\ge 2$. If there exists $(r,s)\in S$ with
$|\theta_{p,q}+\theta_{r,s}|>\delta/2$ then
\begin{eqnarray*}
\frac{\delta}{2} < |\theta_{p,q}+\theta_{r,s}| 
&\le& |\theta_{p,q}|+|\theta_{r,s}|\\
&<& \pi - \delta + \pi\\
&<& 2\pi - \frac{\delta}{2}
\end{eqnarray*}
and we are finished. Otherwise, all $\theta_{r,s}$ with $(r,s)\in S$ satisfy 
$|\theta_{p,q}+\theta_{r,s}|\le\delta/2$; i.e.\ they are all within $\delta/2$ units
of $-\theta_{p,q}$, and so, because $\delta < |\theta_{p,q}| < \pi-\delta$, 
they all have the same sign and satisfy $\delta/2 \le |\theta_{r,s}| \le \pi-\delta/2$.
Now let $(r,s)\in S$ and choose any $(t,u)$ with 
$t\in\{1, 2, \ldots, k\}\setminus\{p,r\}$ and
$u\in\{1, 2, \ldots, k\}\setminus\{q,s\}$. This is possible because $k\ge 3$.
Since $(t,u)\in S$ we have, using the above observations, 
$\delta < |\theta_{r,s}+\theta_{t,u}|<2\pi-\delta$,
which implies the required result.

For the remaining case, we must have 
$|\theta_{p,q}|\in[0,\delta]\cup[\pi-\delta,\pi]$ 
for all $1\le p,q \le k$.  We claim there exist $p$, $q$, $r$, $s$ with
$p\ne r$, $q\ne s$, $|\theta_{p,q}|\in[0,\delta]$, and $|\theta_{r,s}|\in[\pi-\delta,\pi]$. If we prove the claim then we are finished because
\[
\frac{\delta}{2} < \pi-2\delta 
\le \left||\theta_{p,q}|-|\theta_{r,s}|\right| 
\le |\theta_{p,q}+\theta_{r,s}| 
\le |\theta_{p,q}|+|\theta_{r,s}| 
\le \delta+\pi < 2\pi-\frac{\delta}{2}.
\]
Assume for contradiction that the claim is false. By the hypothesis of the
proposition there exist $p$ and $q$ with $|\theta_{p,q}|<\pi-\delta$. Since
$|\theta_{p,q}|\in[0,\delta]\cup[\pi-\delta,\pi]$ for every $1\le p,q\le k$
we must have $|\theta_{p,q}|\in[0,\delta]$. Since we are assuming that the
claim is false, we must have $|\theta_{r,s}|\in[0,\delta]$ for the
$(k-1)^2$ pairs $(r,s)$ with $r\ne p$ and $s\ne q$. But the hypothesis of
the proposition also gives us $(t,u)$ with $|\theta_{t,u}|>\delta$, so an
argument analogous to the previous one shows there must exist $(k-1)^2$ pairs $(v,w)$ with
$|\theta_{v,w}|\in[\pi-\delta,\pi]$. Since $(k-1)^2+(k-1)^2$ exceeds $k^2$, the
total number of ordered pairs, we have a contradiction, as required.
\qed


\proofof{Lemma~\ref{lem:saddle}}
We use the saddlepoint method to estimate the coefficient $C$ in the statement.
Using Cauchy's integral formula, $C$ can
be written in terms of an integral around the product of circles
$z_{p,q}=\rho_{p,q}\exp(i\theta_{p,q})$, $-\pi \le \theta_{p,q}\le \pi$, $(1\le p,q \le k)$, as follows,
\begin{eqnarray}
C &=&\frac{1}{(2\pi i)^{k^2}} \int \frac{     
\exp\left(\frac{1}{2}
\sum_{p\ne r \atop q\ne s}
z_{p,q}z_{r,s}\right)
}{\prod_{p,q} z_{p,q}^{dm_{p,q}n/k+1}}\prod_{p,q}dz_{p,q}\nonumber \\
&=& \frac{1}{(2\pi)^{k^2}\prod_{p,q}\rho_{p,q}^{dm_{p,q}n/k}} \nonumber\\
&&{}\times\int_{\vtheta\in [-\pi,\pi]^{k^2}}
\frac{\exp\left(\frac{1}{2}
\sum_{p\ne r \atop q\ne s}
\rho_{p,q}\rho_{r,s} e^{i(\theta_{p,q}+\theta_{r,s})}\right)
}{\exp(i\sum_{p,q}\theta_{p,q}dm_{p,q}n/k)}
\prod_{p,q}d\theta_{p,q}\lab{eCI}.
\end{eqnarray}
Viewing $\vtheta = \mvec([\theta_{p,q}])$ as a $k^2$-dimensional vector,
let $g(\vtheta)$ denote the integrand in the above expression. Consider
\begin{eqnarray*}
g(\vtheta+\pi\onesvec)
&=&
\frac{\exp\left(\frac{1}{2}
\sum_{p\ne r \atop q\ne s}
\rho_{p,q}\rho_{r,s} e^{i(\theta_{p,q}+\theta_{r,s}+2\pi)}\right)
}{\exp(i\sum_{p,q}\theta_{p,q}dm_{p,q}n/k+ i\pi\sum_{p,q}dm_{p,q}n/k )}\\
&=&
\frac{\exp\left(\frac{1}{2}
\sum_{p\ne r \atop q\ne s}
\rho_{p,q}\rho_{r,s} e^{i(\theta_{p,q}+\theta_{r,s})}\right)
}{\exp(i\sum_{p,q}\theta_{p,q}dm_{p,q}n/k+ i\pi dn )}\\
&=&
g(\vtheta),
\end{eqnarray*}
which holds since $\sum_{p,q}m_{p,q}=k$ and $dn$ is even.
Setting $\delta=\log n / \sqrt{n}$, this tells us that the integral over the
region 
$\{\vtheta \mid |\theta_{p,q}|\le\delta \mbox{ for } 1\le p,q\le k \}$ equals the integral over the region
$\{\vtheta \mid \pi-\delta \le |\theta_{p,q}|\le\pi \mbox{ for } 1\le p,q\le k\}$.
Set each $\rho_{p,q}$ to be the common value
\[
\rho_{p,q} = \rho = \frac{\sqrt{dn}}{k(k-1)}.
\]
We will see that the integral over each of these regions is asymptotic to
\begin{eqnarray}
I &=&
\frac{e^{dn/2}}{\sqrt{2}}
\left(\frac{2\pi}{dn}\right)^{k^2/2}
\frac{ 
k^{k^2}(k-1)^{k(k-1)}
\exp\left( 
\frac{-nd(k-1)^2 
(\mvec \tilde A)^\transpose (J_{k-1}+I_{k-1})^{\otimes 2} \mvec \tilde A
}{2(k^2-2k+2)}
\right) }
{ (k^2-2k+2)^{(k-1)^2/2}
(k-2)^{k-1} } \notag\\
&\ge&
K\exp\left(\frac{dn}{2}-\epsilon'\log^2 n\right)
\label{eq:I}
\end{eqnarray}
(using $a_{p,q} < \epsilon\log n/n^{1/2}$) where $K$ and $\epsilon'$ are constants depending on $\epsilon$.
We will also show that the integral over the remaining region 
is $o(I)$. The lemma then follows by combining these two facts 
with \eqn{eCI}.

To see that the integral over the remaining region is $o(I)$,
let $\vtheta$ be any vector in this region. By the definition of this region
we must have $\min_{p,q}|\theta_{p,q}|<\pi-\delta$ 
and $\max_{p,q}|\theta_{p,q}|> \delta$.
By Lemma~\ref{lem:SaddlePair}
there exist $p^*, q^*, r^*, s^*\in\{1, \ldots, k\}$ 
with $p^*\ne r^*$ and $q^*\ne s^*$ 
such that 
\[
\frac{\delta}{2} \le |\theta_{p^*,q^*}+\theta_{r^*,s^*}| \le 2\pi - \frac{\delta}{2}.
\]
Now
\[
\cos(\theta_{p^*,q^*}+\theta_{r^*,s^*})< \cos\left(\frac{\delta}{2}\right)
= 1 - \frac{\delta^2}{8} + O(\delta^3)
\]
so the absolute value of the integrand is
\begin{eqnarray*}
|g(\vtheta)|
&=&
\frac{\left|\exp\left(\frac{1}{2}
\sum_{p\ne r \atop q\ne s}
\rho^2 e^{i(\theta_{p,q}+\theta_{r,s})}\right)
\right|}{\left|\exp(i\sum_{p,q}\theta_{p,q}dm_{p,q}n/k)\right|}\\
&=&
\exp\left(\frac{1}{2}
\sum_{p\ne r \atop q\ne s}
\rho^2 \cos(\theta_{p,q}+\theta_{r,s})\right)\\
&\le&
\exp\left(\frac{1}{2}
\rho^2 \left((k^2(k-1)^2-1)1 + \cos(\theta_{p^*,q^*}+\theta_{r^*,s^*})\right)\right)\\
&=&
\exp\left(\frac{1}{2}
\rho^2 \left((k^2(k-1)^2-1)1 +  
1 - \frac{\delta^2}{8} + O(\delta^3)\right)
\right)\\
&=&
\exp\left(\frac{1}{2}
\rho^2 k^2(k-1)^2 -  
\rho^2\frac{(\log n)^2}{16n} + O(\rho^2 n^{-3/2}(\log n)^3)
\right)\\
&=&
\exp\left(
\frac{dn}{2}
-\frac{d(\log n)^2}{16k^2(k-1)^2} + o(1)
\right)
\end{eqnarray*}
recalling that we chose $\rho = k^{-1}(k-1)^{-1}\sqrt{dn}$ and
$\delta = \log n / \sqrt{n}$.
Hence $|g(\vtheta)|=o(I)$ if we choose 
$\epsilon$ sufficiently small so that the constant $\epsilon'$ in~\eqref{eq:I} satisfies $\epsilon' < d/\left(16k^2(k-1)^2\right)$.

It remains to show that 
\[
\int_{\vtheta\in [-\delta,\delta]^{k^2}} g(\vtheta) d\vtheta \sim I.
\]
For $\vtheta\in[-\delta,\delta]^{k^2}$ we have
\[
\log g(\vtheta) 
=
\rho^2
\frac{1}{2}
\sum_{p\ne r \atop q\ne s}
\left( 1 + i(\theta_{p,q}+\theta_{r,s}) - \frac{(\theta_{p,q}+\theta_{r,s})^2}{2} + O(|\theta|^3) \right)
- i\frac{dn}{k}\sum_{p,q}\theta_{p,q}m_{p,q}.
\]
Regrouping the terms and substituting $m_{p,q}=k^{-1}+a_{p,q}$ 
this becomes
\begin{eqnarray*}
\log g(\vtheta) 
=&&
\frac{\rho^2}{2}k^2(k-1)^2\\
&+&\sum_{p,q}\theta_{p,q}\left(2i(k-1)^2\frac{\rho^2}{2}-i\frac{dn}{k}\left(\frac{1}{k}+a_{p,q}\right)\right) \\
&-&\frac{\rho^2}{2}\left((k-1)^2\sum_{p,q}\theta_{p,q}^2 + 
\sum_{p\ne r \atop q\ne s}\theta_{p,q}\theta_{r,s}\right)\\
&+& O(\rho^2|\vtheta|^3)
\end{eqnarray*}
Let $c$ be the constant 
$c=d/(2k^2(k-1)^2)$.
Recalling $\rho = k^{-1}(k-1)^{-1}\sqrt{dn}$ 
we find 
\bel{elgg}
\log g(\vtheta) 
=
\frac{dn}{2}
-i\frac{dn}{k}(\mvec A)^\transpose \vtheta
-cn\vtheta^\transpose B \vtheta
+ O(n^{-1/2}(\log n)^3)
\ee
where $B$ is the matrix
\[
B = (k-1)^2 I_{k^2} + (J_k-I_k)^{\otimes 2}.
\]
Define $h(\theta) = -i(dn/k)(\mvec A)^\transpose\vtheta-cn\vtheta^\transpose B\vtheta$.
Lemma~\ref{lem:Evec2} gives us an orthonormal basis $\{f^{(p,q)}\}_{p,q=1}^k$ 
of eigenvectors for $B$ and corresponding sequence of
eigenvalues $(\lambda_{p,q})_{p,q=1}^k$.
 Introduce the new variables $(\tau_{p,q})_{p,q=1}^k$
to perform the change of basis $\vtheta=\sum_{p,q}f^{(p,q)}\tau_{p,q}$. This
gives
\begin{eqnarray*}
h(\vtheta)
&=&
-i\frac{dn}{k}(\mvec A)^\transpose\sum_{p,q}f^{(p,q)} \tau_{p,q}
-cn\sum_{p,q}\lambda_{p,q}\tau_{p,q}^2\\
&=&
\sum_{p,q}\left(
-i(dn/k)(\mvec A)^\transpose f^{(p,q)} \tau_{p,q}
-cn\lambda_{p,q}\tau_{p,q}^2
\right).
\end{eqnarray*}
Let $p,q\in\{1,\ldots,k\}$. 
Using the identity
\[
\int_{-\infty}^\infty e^{ax-bx^2}dx = \sqrt{\frac{\pi}{b}}\exp\left(\frac{a^2}{4b}\right)
\]
(for $b>0$), we have that $\int_{[-\infty,\infty]^{k^2}}
\exp(h(\vtheta)) d\vtheta$ is
a product of terms of the form
\begin{eqnarray*}
&&
\int_{-\infty}^{\infty}\exp\left(
-i\frac{dn}{k}(\mvec A)^\transpose f^{(p,q)} \tau_{p,q}
-cn\lambda_{p,q}\tau_{p,q}^2
\right)d\tau_{p,q}\\
&=&
\sqrt{\frac{\pi}{cn\lambda_{p,q}}}
\exp\left( \frac{-d^2n((\mvec A)^\transpose f^{(p,q)})^2 }
{4ck^2\lambda_{p,q}}  \right).
\end{eqnarray*}
So by Lemma~\ref{lGauss}, for some constant $c'>0$ we have
\begin{eqnarray*}
\int_{[-\delta,\delta]^{k^2}} 
\exp(h(\vtheta))d\vtheta
&=&
\int_{[-\infty,\infty]^{k^2}} 
\exp(h(\vtheta))d\vtheta
+O(e^{-c'(\log n)^2})\\
&=&
\prod_{p,q}
\sqrt{\frac{\pi}{cn\lambda_{p,q}}}
\exp\left( \frac{-d^2n((\mvec A)^\transpose f^{(p,q)})^2 }
{4ck^2\lambda_{p,q}}  \right)
+O(e^{-c'(\log n)^2})\\
&\sim&
\prod_{p,q}
\sqrt{\frac{\pi}{cn\lambda_{p,q}}}
\exp\left( \frac{-d^2n((\mvec A)^\transpose f^{(p,q)})^2 }
{4ck^2\lambda_{p,q}}  \right)
\end{eqnarray*}
since the entries of $A$ are less than $\epsilon n^{-1/2}\log n$, $c'$ does not depend on $\epsilon$ and we will
choose $\epsilon$ to be sufficiently small.
 Recalling \eqn{elgg} we now have
\bel{eIntProd}
\int_{[-\delta,\delta]^{k^2}} g(\vtheta) d\vtheta
\sim
e^{dn/2}
\prod_{p,q}
\sqrt{\frac{\pi}{cn\lambda_{p,q}}}
\exp\left( \frac{-d^2n((\mvec A)^\transpose f^{(p,q)})^2 }
{4ck^2\lambda_{p,q}}  \right)
\ee
We will simplify the above product using the values of $\lambda_{p,q}$ given in Lemma~\ref{lem:Evec2}.
First, the contribution to  the product from  $1\le p,q \le k-1$ is
\begin{eqnarray*}
&&\prod_{p=1}^{k-1}\prod_{q=1}^{k-1}
\sqrt{\frac{\pi}{cn\lambda_{p,q}}}
\exp\left( \frac{-d^2n((\mvec A)^\transpose f^{(p,q)})^2 }
{4ck^2\lambda_{p,q}}  \right)\\
&\sim&
\left(\sqrt{\frac{\pi}{cn(k^2-2k+2)}}\right)^{(k-1)^2}
\prod_{p=1}^{k-1}\prod_{q=1}^{k-1}
\exp\left( \frac{-d^2n((\mvec A)^\transpose f^{(p,q)})^2 }
{4ck^2(k^2-2k+2)}  \right)\\
&=&
\left(\sqrt{\frac{\pi}{cn(k^2-2k+2)}}\right)^{(k-1)^2}
\exp\left( \frac{-d^2n( \sum_{p=1}^{k-1}\sum_{q=1}^{k-1}
(\mvec A)^\transpose f^{(p,q)})^2 }
{4ck^2(k^2-2k+2)}  \right)\\
&=&
\left(\sqrt{\frac{\pi}{cn(k^2-2k+2)}}\right)^{(k-1)^2}
\exp\left( \frac{-d^2n 
(\mvec \tilde A)^\transpose (J_{k-1}+I_{k-1})^{\otimes 2} \mvec \tilde A
}
{4ck^2(k^2-2k+2)}  \right)\\
\end{eqnarray*}
where the last step used Lemmata \ref{lem:Evec3}(b) and \ref{lem:Evec3}(c). 
The contribution to the product when exactly one of $p$ or $q$ equals $k$
is
\[
\left(\sqrt{\frac{\pi}{cn(k-1)(k-2)}}\right)^{2(k-1)}
\]
since Lemma~\ref{lem:Evec3}(a) tells us that $((\mvec A)^\transpose f^{(p,q)})^2=0$ when $p=k$ or $q=k$. When $p=q=k$ the contribution to the product is
\[
\sqrt{\frac{\pi}{2cn(k-1)^2}}
\]
Substituting these contributions into \eqn{eIntProd} we get
\begin{eqnarray*}
&&\int_{[-\delta,\delta]^{k^2}} g(\vtheta) d\vtheta\\
&\sim&
e^{dn/2}
\left(\frac{\pi}{cn}\right)^{k^2/2}
\frac{ \exp\left( 
\frac{-nd^2 
(\mvec \tilde A)^\transpose (J_{k-1}+I_{k-1})^{\otimes 2} \mvec \tilde A
}{4ck^2(k^2-2k+2)}
\right) }
{ (k^2-2k+2)^{(k-1)^2/2}
\left((k-1)(k-2)\right)^{k-1}
\sqrt{2}(k-1)}\\
&=&
e^{dn/2}
\left(\frac{2\pi k^2 (k-1)^2}{dn}\right)^{k^2/2}
\frac{ \exp\left( 
\frac{-nd(k-1)^2 
(\mvec \tilde A)^\transpose (J_{k-1}+I_{k-1})^{\otimes 2} \mvec \tilde A
}{2(k^2-2k+2)}
\right) }
{ (k^2-2k+2)^{(k-1)^2/2}
\left((k-1)(k-2)\right)^{k-1}
\sqrt{2}(k-1)}\\
&=&
\frac{e^{dn/2}}{\sqrt{2}}
\left(\frac{2\pi}{dn}\right)^{k^2/2}
\frac{ 
k^{k^2}(k-1)^{k(k-1)}
\exp\left( 
\frac{-nd(k-1)^2 
(\mvec \tilde A)^\transpose (J_{k-1}+I_{k-1})^{\otimes 2} \mvec \tilde A
}{2(k^2-2k+2)}
\right) }
{ (k^2-2k+2)^{(k-1)^2/2}
(k-2)^{k-1} }\\
&=& I,
\end{eqnarray*}
as required.
\qed

\section{\ldots and for $\mbf n$ not divisible by $\mbf k$}
\label{sec:general_n}
Define $k'=2k$ if $dk$ is odd or $k'=k$ otherwise. Note that the conditions we assumed so far (i.e.~$n$ divisible by $k$ and $dn$ even) can be rewritten as simply $n\equiv 0$ (mod $k'$). Therefore we only need to consider the case $n\equiv r$ (mod $k'$) for each integer $r$ such that $0<r<k'$ and $dr$ is even.
One possibility is to rework the whole argument of this paper but with slightly unbalanced colourings. Instead, the asymmetry in the argument can be somewhat reduced by using an argument relating different models of random regular graphs. We first treat the case $n\equiv 0$ (mod $k'$) in more depth, and prove the following.
\begin{thm}
\lab{TMAIN2}
Fix nonnegative integers $d\ge3$, $k$ and $\ell$ such that $d<2(k-1)\log(k-1)$. Consider the $d$-regular graphs with $n$ vertices ($n$ divisible by $k$ and $dn$ even) and a distinguished ordered set of $\ell$ edges whose endpoints induce a perfect matching (i.e.~no two edges in the distinguished set are adjacent to the same edge or incident with the same vertex). Let $G$ be chosen uniformly at   random from such structures. Then $G$ a.a.s.\ has a $k$-colouring in which all $\ell$ distinguished edges have end vertices coloured $1$ and $2$.
\end{thm}
\proof
Consider the probability space $\Omega_{n,d,\ell}$ with uniform probability distribution, and whose underlying set consists of pairings in $\Pnd$ with an ordered set $L$ of $\ell$ distinguished pairs of points, such that no two pairs in $L$ are incident with the same vertex. Let $\hat Y$ denote the  number of balanced $k$-colourings of a pairing containing $L$, in which the distinguished pairs join vertices of colours $1$ and $2$. We will show that
\bel{yhat}
\ex \hat Y \sim \binom{k}{2}^{-\ell}\ex Y,
\ee
that~\eqref{E3} holds with $Y$ replaced by $\hat Y$ (and no other adjustment), and that
\bel{yhatsq}
\ex (\hat Y^2) \sim \binom{k}{2}^{-2\ell}\ex (Y^2).
\ee
The theorem  then follows immediately by the argument in the last few sentences of the proof of Theorem~\ref{thm:main} (for $n$ divisible by $k$).

To show~\eqref{yhat} and the analogue of~\eqref{E3}, we define for integers $r\ge0$ and $p_1,\ldots,p_r\ge0$
$N(p_1,\ldots,p_r)$ to be the set of triples $(P,C,\Gamma)$ such that $P$ is a pairing in $\Pnd$, $C$ is a balanced $k$-colouring of $P$ and $\Gamma=(\Gamma_{i,j})_{1\le i\le r\atop 1\le j\le p_i}$ is an ordered $(p_1+\cdots+p_r)$-tuple of different cycles in $P$. In view of that, we can express
\begin{gather}
\ex(Y[X_1]_{p_1}\cdots [X_r]_{p_r}) = \frac{|N(p_1,\ldots,p_r)|}{|\Pnd|}
\quad\text{and}
\notag\\
\ex(\hat Y[X_1]_{p_1}\cdots [X_r]_{p_r}) = \frac{1}{|\Omega_{n,d,\ell}|}
\sum_{(P,C,\Gamma)\in N(p_1,\ldots,p_r)} h(P,C),
\label{eq:jointhat}
\end{gather}
where $h(P,C)$ gives the number of choices of the ordered set $L$ of $\ell$ pairs in $P$ that have the required colours at their ends.
Almost all triples in $N(p_1,\ldots,p_r)$ correspond to pairings with $dn/(k(k-1)) \pm O(n^{1/2} \log n)$ edges between each two colour classes, while the remaining triples contribute with at most a $O(n^{-\Theta(\log n)})$ factor of the total.
To verify this claim, observe that for each fixed $C$ and $\Gamma$, the number of pairings $P\in\Pnd$ compatible with $C$ and $\Gamma$ has a factor $F=\sum_{\{b_{i,j}\}}1/(\prod_{1\le i<j\le k}b_{i,j}!)$, where $b_{i,j}$ denotes the number of edges between colour classes $i$ and $j$ excluding the edges of the cycles in $\Gamma$ (see the computation of $\ex Y$ and $\ex(YX_1)$ in Section~\ref{sec:joint}).
After using Stirling's formula to estimate the factorials in $F$, it is easy to check that the main contribution to $F$ is due to terms with all $b_{i,j} = dn/(k(k-1)) \pm O(n^{1/2} \log n)$ and that the weight of the remaining terms is $O(F / n^{\Theta(\log n)})$ as required.
Next, observe that $h(P,C)\sim (dn/(k(k-1)))^\ell$ if $P$ has $dn/(k(k-1)) \pm O(n^{1/2} \log n)$ edges with endpoints coloured $1$ and $2$, and that $h(P,C)$ is always $O(n^\ell)$. Therefore, we can estimate the sum in the right side of~\eqref{eq:jointhat} and combine it with $|\Omega_{n,d,\ell}|\sim|\Pnd|(dn/2)^\ell$ to deduce
\[
\ex(\hat Y[X_1]_{p_1}\cdots [X_r]_{p_r}) \sim \binom{k}{2}^{-\ell}
\ex(Y[X_1]_{p_1}\cdots [X_r]_{p_r}),
\]
as required for~\eqref{yhat} and the analogue of~\eqref{E3}.

Similarly, to estimate $\ex (\hat Y^2)$ we define $N'$ to be the set of triples $(P,C_1,C_2)$ such that $P$ is a pairing in $\Pnd$ and $C_1$, $C_2$ are balanced $k$-colourings of $P$, and write
\bel{eq:secondhat}
\ex(Y^2) = \frac{|N'|}{|\Pnd|}
\qquad\text{and}\qquad
\ex(\hat Y^2) = \frac{1}{|\Omega_{n,d,\ell}|}
\sum_{(P,C_1,C_2)\in N'} h'(P,C_1,C_2),
\ee
where $h'(P,C_1,C_2)$ gives the number of choices of the ordered set $L$ of $\ell$ pairs in $P$ that have the required colours at their ends in both colourings $C_1$ and $C_2$. Given $P\in\Pnd$, let $b_{pqrs}$ denote the number of edges between points coloured $(p,q)$ and points coloured $(r,s)$. We will show that almost all triples in $N'$ correspond to pairings with $b_{pqrs}=dn/(k^2(k-1)^2)\pm O(n^{1/2}\log n)$, and that the remaining triples are at most a $O(n^{-\Theta(\log n)})$ fraction of the total. Then~\eqref{yhatsq} follows immediately from~\eqref{eq:secondhat}, since $h'(P,C_1,C_2)$ is always $O(n^\ell)$ and also $h'(P,C_1,C_2) \sim (2dn/(k^2(k-1)^2))^\ell$ whenever $b_{1122}\sim b_{1221}\sim dn/(k^2(k-1)^2)$.
To prove the remaining claim, we first recall from the last lines of the proof of Proposition~\ref{prop:moments}$(b)$ that we can restrict our attention to triples $(P,C_1,C_2)$ with colour count $M=[m_{p,q}]$ where $|m_{p,q}-1/k|<\epsilon \log n/n^{1/2}$, since all other triples contribute $O(|N'|/n^{\Theta(\log n)})$ to $|N'|$.
Observe that while counting $|N'|$ we encounter a factor $F'=\sum_{\cB(M,n)} \prod_{1\le p,q,r,s\le k\atop p<r,q\ne s} 1/b_{pqrs}!$. We can easily bound the weight in $F'$ due to terms in which some $b_{pqrs}$ is not $dn/(k^2(k-1)^2)\pm O(n^{1/2}\log n)$, and find an extra factor of $O(n^{-\Theta(\log n)})$ compared to the estimation of $F'$ in Lemma~\ref{lem:saddle}.
\qed

\proofof{Theorem~\ref{thm:main} (for $\mbf n$ not divisible by $\mbf k$)}
It only remains to show that if $d\ge3$ and $d<2(k-1)\log(k-1)$, then $\Gnd$ is a.a.s.\ $k$-colourable for $n$ not divisible by $k$. We use the type of argument employed at the end of Section~3 of~\cite{RW94}.

Recall the definition of $k'$ in the beginning of the section, and let $r$ be any integer such that $0<r<k'$ and $dr$ is even. Suppose $n\equiv r$ (mod $k'$). Take a random $d$-regular graph $G$ with $n$ vertices,  and assume that the first $r$ vertices $v_1,v_2,\ldots,v_r$ are at distance at least $4$. This happens a.a.s.\ and we simply discard $G$ if this property fails to hold. Delete $v_1,v_2,\ldots,v_r$, and join up the former $dr$ neighbours of these vertices by a random perfect matching, which is added to what remains of graph $G$. Leave the $dr/2$ added edges as a distinguished ordered set of edges (in any random order), and observe that none of these edges are adjacent to each other by construction. It is easy to show and well known that a given vertex of a random $d$-regular graph is a.a.s.\ not in a cycle of length less than 4 (or 100, for that matter). It follows that a.a.s.\ no multiple edges occur due to the new edges.  Throw the graph away if this last property fails to hold. The result is a random $d$-regular graph  on $n-r\equiv 0$ (mod $k'$) vertices with an ordered set of $dr/2$ distinguished edges, no two adjacent to the same edge or incident with the same vertex. Let us call this $G'$.

The distribution of $G'$ is uniform, since for each $G'$ the number of ways of reinstating the edges to $v_1,v_2,\ldots,v_r$ is exactly $\binom{dr}{d,d,\ldots,d}$.
Thus, by Theorem~\ref{TMAIN2}, $G'$ a.a.s.\ has a $k$-colouring such that the $dr/2$ distinguished edges join vertices of colours $1$ and $2$. We can use exactly this colouring on $V(G)\setminus\{v_1,v_2,\ldots,v_r\}$, and colour $v_1,v_2,\ldots,v_r$ with colour $3$ to obtain a $k$-colouring of $G$. (Note that $k\ge3$ from our assumptions on $d$ and $k$). \qed

\bibliographystyle{abbrv}
\bibliography{references_short}

\end{document}